\RequirePackage{fix-cm} 
\documentclass[a4paper, twoside, reqno, dvips, 12pt]{amsart}
\usepackage{fixltx2e}   

\usepackage[latin1]{inputenc}
\usepackage[T1]{fontenc}

\usepackage{amsbsy}
\usepackage{dsfont}
\usepackage{xspace}
\usepackage{amsgen}
\usepackage{amsthm}
\usepackage{amssymb}
\usepackage{amsmath}
\usepackage{upgreek}
\usepackage{amsfonts}
\usepackage{extarrows}
\usepackage{mathtools}
\usepackage{textgreek}
\usepackage[nice]{nicefrac}

\usepackage{mathabx}
\usepackage{relsize}
\usepackage{textcomp}
\usepackage[mathcal]{euscript}

\usepackage{mathrsfs}
\DeclareMathAlphabet{\mathscrbf}{OMS}{mdugm}{b}{n}

\usepackage[scr]{rsfso}

\usepackage{a4wide}

\usepackage[dvipsnames, table]{xcolor}

\definecolor{bckg}{RGB}{20.8, 20.8, 20.8}
\definecolor{oneblue}{rgb}{0.0, 0.0, 0.85}
\definecolor{Lightblue}{RGB}{214, 214, 214}
\definecolor{bluepigment}{rgb}{0.2, 0.2, 0.6}
\definecolor{charcoal}{rgb}{0.21, 0.27, 0.31}
\definecolor{denimblue}{rgb}{0.08, 0.38, 0.74}
\definecolor{Lightgray}{rgb}{0.89, 0.89, 0.89}
\definecolor{darkgrey}{rgb}{0.273, 0.281, 0.30}
\definecolor{darkelectricblue}{rgb}{0.33, 0.41, 0.47}

\usepackage[sort&compress, comma, square, numbers]{natbib}

\usepackage{psfrag}
\usepackage{graphicx}
\usepackage{subfigure}
\usepackage{morefloats}
\usepackage{indentfirst}

\usepackage{acronym}
\usepackage{microtype}
\usepackage[labelsep=period,%
            labelfont={bf,sf,color=bluepigment},%
            justification=raggedright]{caption}

\usepackage[perpage, symbol]{footmisc}

\usepackage[usenames, dvipsnames, pdf]{pstricks}
\usepackage{epsfig}
\usepackage{pst-grad} 
\usepackage{pst-plot} 

\usepackage{multirow}
\usepackage{tabularx}

\usepackage[colorlinks,
           urlcolor=oneblue,
           linkcolor=denimblue,
           citecolor=NavyBlue,
           bookmarksopen=false,
           pdfpagemode=UseNone,
           pagebackref]{hyperref}

\usepackage[explicit]{titlesec}

\titleformat{\section}
  {\color{NavyBlue}\Large\sffamily\bfseries}
  {}
  {0em}
  {\colorbox{bckg!5}{\parbox{\dimexpr\linewidth-2\fboxsep\relax}{\centering\thesection. #1}}}
  [\vspace*{0.33em}]

\titleformat{name=\section,numberless}
  {\color{NavyBlue}\Large\sffamily\bfseries}
  {}
  {0.0em}
  {\colorbox{bckg!5}{\parbox{\dimexpr\linewidth-2\fboxsep\relax}{\centering#1}}}
  [\vspace*{0.33em}]

\titleformat{\subsection}
  {\color{NavyBlue}\large\sffamily\bfseries}
  {}
  {0.0em}
  {\colorbox{bckg!5}{\parbox{\dimexpr\linewidth-2\fboxsep\relax}{\centering\thesubsection. #1}}}
  [\vspace*{0.33em}]

\titleformat{name=\subsection,numberless}
  {\color{NavyBlue}\Large\sffamily\bfseries}
  {}
  {0em}
  {\colorbox{bckg!5}{\parbox{\dimexpr\linewidth-2\fboxsep\relax}{\centering#1}}}
  [\vspace*{0.33em}]

\titleformat{\subsubsection}
  {\color{bluepigment}\sffamily\normalsize\bfseries}
  {\thesubsubsection}
  {0.5em}
  {\colorbox{bckg!5}{\parbox{\dimexpr\linewidth-16\fboxsep\relax}{\centering#1}}}
  [\vspace*{0.33em}]

\titleformat{\paragraph}[runin]
  {\color{bluepigment}\sffamily\small\bfseries}
  {}
  {0em}
  {#1}

\titlespacing{\section}{1.0em}{1.5em plus 2pt minus 2pt}%
{1.0em plus 2pt minus 2pt}[0em]
\titlespacing{\subsection}{1.0em}{1.5em plus 2pt minus 2pt}%
{1.0em}[0em]
\titlespacing{\subsubsection}{1.0em}{1.5em plus 2pt minus 2pt}%
{1.0em plus 2pt minus 2pt}[0em]

\usepackage{titletoc}

\contentsmargin{0.5em}
\newlength{\tocsep} 
\titlecontents{section}[\tocsep]
  {\addvspace{10pt}\bfseries\sffamily}
  {\contentslabel[\thecontentslabel]{\tocsep}}
  {}
  {\ \titlerule*[0.75pc]{.}\ \thecontentspage}
  []
\titlecontents{subsection}[\tocsep]
  {\addvspace{8pt}\sffamily}
  {\contentslabel[\thecontentslabel]{\tocsep}}
  {}
  {\ \titlerule*[0.5pc]{.}\ \thecontentspage}
  []
\titlecontents*{subsubsection}[\tocsep]
  {\addvspace{2pt}\footnotesize\sffamily}
  {}
  {}
  {\ \titlerule*[0.35pc]{.}\ \thecontentspage}
  [\\*]

\makeatletter
\def\@setauthors{%
  \begingroup
  \def\thanks{\protect\thanks@warning}%
  \trivlist
  \centering\footnotesize \@topsep30\p@\relax
  \advance\@topsep by -\baselineskip
  \item\relax
  \author@andify\authors
  \def\\{\protect\linebreak}%
  \textsc{\normalsize\textcolor{darkelectricblue}{\authors}}%
  \ifx\@empty\contribs
  \else
    ,\penalty-3 \space \@setcontribs
    \@closetoccontribs
  \fi
  \endtrivlist
  \endgroup
}
\def\@settitle{\begin{center}%
  \baselineskip14\p@\relax
    \bfseries
    \textsc{\Large\textcolor{charcoal}{\@title}}
  \end{center}%
}
\makeatother

\usepackage{enumitem}
\setlist[description]{%
  topsep=30pt,               
  itemsep=5pt,               
  font={\bfseries\sffamily\color{NavyBlue}}, 
}

\usepackage{fancyhdr}
\usepackage{lastpage}

\newcommand*\Title{\textcolor{bluepigment}{Analysis and improvement of the VTT model}}
\newcommand*\Authors{\textcolor{bluepigment}{J.~Berger, H.~Le Meur, D.~Dutykh, \etal}}
\newcommand*{\plogo}{\textcolor{gray}{{\texttt{arXiv.org} / \textsc{hal}}}} 

\numberwithin{equation}{section}

\newcommand{\ie}{\emph{i.e.}\xspace}

\newcommand{\etal}{\emph{et al.}\xspace}

\renewcommand{\sim}{\thicksim}

\newcommand{\pd}[2]{\frac{\partial #1}{\partial\/ #2}}
\newcommand{\od}[2]{\frac{\mathrm{d} #1}{\mathrm{d}\/#2}}

\newcommand{\eqdef}{\mathop{\stackrel{\,\mathrm{def}}{:=}\,}}

\makeatletter
\renewcommand*\env@matrix[1][\arraystretch]{%
  \edef\arraystretch{#1}%
  \hskip -\arraycolsep
  \let\@ifnextchar\new@ifnextchar
  \array{*\c@MaxMatrixCols c}}
\makeatother

\newcommand*\egal{\ = \ }
\newcommand*\plus{\ + \ }
\newcommand*\moins{\ - \ }

\newcommand{\Mmax}{M_{\,\mathrm{max}}}
\newcommand{\Minf}{M_{\,\infty}}
\newcommand{\Mzero}{M_{\,0}}
\newcommand{\phic}{\phi_{\,\mathrm{c}}}
\newcommand{\gP}{\mathbf{P}}
\newcommand{\gM}{\mathbf{M}}
\newcommand{\gF}{\mathds{F}}
\renewcommand{\O}{\mathcal{O}}

\newcommand{\apriori}{\emph{a priori}$\;$}

\begin{document}

\title[\Title]{Analysis and improvement of the VTT mold growth model: application to bamboo fiberboard}

\author[J.~Berger]{Julien Berger$^*$}
\address{\textbf{J.~Berger:} Univ. Grenoble Alpes, Univ. Savoie Mont Blanc, UMR 5271 CNRS, LOCIE, 73000 Chamb\'ery, France}
\email{Julien.Berger@puniv-smb.fr}
\urladdr{https://www.researchgate.net/profile/Julien\_Berger3/}
\thanks{$^*$ Corresponding author}

\author[H.~Le Meur]{Herv\'e Le Meur}
\address{\textbf{H.~Le Meur:} Laboratoire de Math\'ematiques d'Orsay, Univ. Paris-Sud, CNRS, Universit\'e Paris-Saclay, 91405 Orsay, France}
\email{Herve.LeMeur@math.u-psud.fr}
\urladdr{https://www.math.u-psud.fr/~lemeur/}

\author[D.~Dutykh]{Denys Dutykh}
\address{\textbf{D.~Dutykh:} Univ. Grenoble Alpes, Univ. Savoie Mont Blanc, CNRS, LAMA, 73000 Chamb\'ery, France and LAMA, UMR 5127 CNRS, Universit\'e Savoie Mont Blanc, Campus Scientifique, 73376 Le Bourget-du-Lac Cedex, France}
\email{Denys.Dutykh@univ-smb.fr}
\urladdr{http://www.denys-dutykh.com/}

\author[D.~M.~Nguyen]{Dang Mao Nguyen}
\address{\textbf{D.~M.~Nguyen:} Univ. Grenoble Alpes, Univ. Savoie Mont Blanc, UMR 5271 CNRS, LOCIE, 73000 Chamb\'ery, France}

\author[A.-C.~Grillet]{Anne-C\'ecile Grillet}
\address{\textbf{A.-C.~Grillet:} Univ. Grenoble Alpes, Univ. Savoie Mont Blanc, UMR 5271 CNRS, LOCIE, 73000 Chamb\'ery, France}
\email{Anne-Cecile.grillet@univ-smb.fr}



\begin{titlepage}
\thispagestyle{empty} 
\noindent
{\Large Julien \textsc{Berger}}\\
{\it\textcolor{gray}{Polytech Annecy--Chamb\'ery, LOCIE, France}}
\\[0.02\textheight]
{\Large Herv\'e \textsc{Le Meur}}\\
{\it\textcolor{gray}{Universit\'e Paris-Saclay, Orsay, France}}
\\[0.02\textheight]
{\Large Denys \textsc{Dutykh}}\\
{\it\textcolor{gray}{CNRS, Universit\'e Savoie Mont Blanc, France}}
\\[0.02\textheight]
{\Large Dang Mao \textsc{Nguyen}}\\
{\it\textcolor{gray}{Polytech Annecy--Chamb\'ery, LOCIE, France}}
\\[0.02\textheight]
{\Large Anne-C\'ecile \textsc{Grillet}}\\
{\it\textcolor{gray}{Polytech Annecy--Chamb\'ery, LOCIE, France}}
\\[0.05\textheight]

\colorbox{Lightblue}{
  \parbox[t]{1.0\textwidth}{
    \centering\huge\sc
    \vspace*{0.75cm}
    
    \textcolor{bluepigment}{Analysis and improvement of the VTT mold growth model: application to bamboo fiberboard}
    
    \vspace*{0.75cm}
  }
}

\vfill 

\raggedleft     
{\large \plogo} 
\end{titlepage}


\newpage
\thispagestyle{empty} 
\par\vspace*{\fill}   
\begin{flushright} 
{\textcolor{denimblue}{\textsc{Last modified:}} \today}
\end{flushright}


\newpage
\maketitle
\thispagestyle{empty}


\begin{abstract}

The reliability of a model is its accuracy in predicting the physical phenomena using the known input parameters. It also depends on the model's ability to estimate relevant parameters using observations of the physical phenomena. In this paper, the reliability of the VTT model is investigated under these two criteria for various given temperature and relative humidity constant in time. First of all, experiments are conducted on bamboo fiberboard. Using these data, five parameters of the VTT model, defining the mold vulnerability class of a material, are identified. The results highlight that the determined parameters are not within the range of the classes defined in the VTT model. In addition, the quality of the parameter estimation is not satisfactory. Then the sensitivity of the numerical results of the VTT model is analyzed by varying an input parameter. These investigations show that the VTT mathematical formulation of the physical model of mold growth is not reliable. An improved model is proposed with a new mathematical formulation. It is inspired by the logistic equation whose parameters are estimated using the experimental data obtained. The parameter estimation is very satisfactory. In the last parts of the paper, the numerical predictions of the improved model are compared to experimental data from the literature to prove its reliability.

\bigskip
\noindent \textbf{\keywordsname:} Mold growth prediction models; VTT model; Model reliability; Parameter estimation problem; Structural identifiability; Fisher matrix \\


\end{abstract}


\newpage
\tableofcontents
\thispagestyle{empty}


\newpage
\section{Introduction}

\subsection{Context}

Excessive humidity in buildings damages the construction quality and affects indoor air quality as well as the occupants' thermal comfort. Moisture is a source of several disorders in buildings, as reported in \cite{Berger2015a, Harris2001}. For example, mold growth is a consequence of high moisture levels and is of capital importance since it can be toxic for the occupants of the building, causing allergies, diseases or infections \cite{Palot2011, Reboux2010, Freire2017}. It also causes the development of other moisture-related damages such as metal corrosion or degradation of the materials by chemical reactions.

Several models of mold development exist and have been listed by \textsc{Vereecken} in \cite{Vereecken2015}. They can be conventionally divided into two classes: the static and the dynamic models. The former only indicates the initiation of the biological process while the latter represents the dynamic of the physical phenomena of mold growth and decline. The so-called VTT model is dynamic and one of the most used \cite{Vereecken2015}. Interested readers are invited to consult \cite{Fedorik2013, Harrestrup2016} for examples of its practical applications. This model depends on the time-dependent relative humidity and temperature conditions that can be calculated from a building Heat And Moisture (HAM) simulation program \cite{Woloszyn2008, Mendes2008}.


\subsection{Problem statement}
\label{sec:problem_statement}

For economic and social reasons, it is essential to have robust models able to represent the physical phenomena of mold growth on building materials. In \cite{Vereecken2012}, \textsc{Vereecken} and \textsc{Roels} highlighted a number of discrepancies between the numerical results of the VTT model and experimental data measured in \cite{Johansson2013}. Similar conclusions have been drawn in recent studies. Most particularly, in \cite{Colinart2017}, the reliability of numerical predictions obtained for cyclic conditions of relative humidity and temperature is questioned. In \cite{Marincioni2017}, the suitability of the mold growth model is investigated for prediction on woodfiber insulation on an internal wall using measured temperature and relative humidity data. Even if improvements have been made, as mentioned in \cite{Viitanen2015}, based on these results, it is of major importance to develop more reliable models.

The robustness of a model is also based on the possibility of accurately estimating the relevant parameters for new innovative materials. These parameters are incorporated in a so-called \emph{vulnerability}\footnote{The VTT model defines four classes of material \emph{sensitivity} to mold growth. However, in order to avoid any confusion with the mathematical analysis performed in this study, their word \emph{sensitivity} is replaced by \emph{vulnerability} here. In this study, the model \emph{sensitivity} defines the variation of the model prediction due to changes of its parameters.} class in the VTT model. This paper aims at analyzing the robustness of the VTT model according to both aspects: (i) the estimation of parameters using experimental observations and (ii) the prediction of the physical phenomena. These investigations are carried for temperature and relative humidity kept constant in time. Even if the VTT model has been proposed for time-varying temperature and relative humidity, its reliability should be ensued for at least constant conditions.

The manuscript is organized as follows. Section~\ref{sec:modelVTT} reviews the formulation of the VTT mold growth model. Both the identifiability and the identification of the parameters are presented together with mathematical tools to quantify the potential error in parameter identification. In Section~\ref{sec:experimental_data}, the experimental data obtained for bamboo fiber are presented. It appears that they may not fit the mathematical model driving us to project these observed data. The numerical results of the parameter estimation for the VTT model are presented in Section~\ref{sec:estimation_parameters_VTT} with a quantification of the potential confidence. Section~\ref{sec:Reliability_VTT_direct} evaluates the VTT model's ability to predict mold growth. The precision of the numerical predictions is analyzed according to small changes in the parameters' numerical values. This discussion will be followed by comments on the robustness of the mathematical formulation of the VTT model, with a proposal of a better mathematical model, in Section~\ref{sec:improvement_model}, with better numerical properties. Last, we draw the main conclusions of this study in Section~\ref{sec:final_remarks}.


\section{Methodology}
\label{sec:modelVTT}

\subsection{Mathematical formulation of the VTT mold growth model}

A mathematical formulation of the physical phenomena for wood-based material was first proposed by \textsc{Hukka} and \textsc{Viitanen} in \cite{Hukka1999}. Its extension to other building materials was suggested by \textsc{Viitanen} \cite{Viitanen2010}. The field of interest is the quantity $M \ \in \ \bigl[\, 0 \,,\, 6\,\bigr]\,$, describing the amount of mold at the surface of the material and computed by solving the initial value (or \textsc{Cauchy}) problem for the following differential equation: 
\begin{align}
\label{eq:odeM}
& \frac{\mathrm{d}M}{\mathrm{d}t} \egal \frac{k_{\,1}\,\bigl(\,M(\,t\,)\,\bigr) \cdot k_{\,2} \,\bigl(\,M(\,t\,)\,\bigr)}{ f \, \bigl(\, \,T \,(\,t\,) \,,\, \phi \,(\,t\,)\,\bigr)} \,, && t \ > \ 0 \,,
\end{align}
where $f\, (\,T \,,\, \phi \,)$ is a known function of the temperature $T$ (in $\mathsf{^{\,\circ} \, Celsius}$), the relative humidity $\phi$ (in percentage) and other parameters such as the material's surface quality. Time is expressed in hours ($\mathsf{h}$). The initial condition is $M(\,t\egal0\,) \egal 0\,$. The quantity $M$ is dimensionless. For a surface quality corresponding to a sawn surface, the function $f$ can be written as: 
\begin{align}
\label{eq:function_f}
f\, (\,T \,,\, \phi \,)
&  \egal b_{\,0} \, \exp \biggl(\, b_{\,1} \, \ln \bigr[\, T \,(\,t\,) \,\bigl] 
\plus b_{\,2} \, \ln \bigr[\, \phi \,(\,t\,) \,\bigl] 
\plus b_{\,3}  \,\biggr) \,, \\
 b_{\,0} & \egal 168 \,, 
\qquad \ \qquad b_{\,1} \egal -0.68 \,,\nonumber \\
b_{\,2} & \egal -13.9 \,, 
\qquad \ \quad b_{\,3} \egal 66.02 \,. \nonumber
\end{align}

The coefficients $k_{\,1}$ and $k_{\,2}$ are defined as: 
\begin{align*}
k_{\,1}\,(\,M\,) &\egal 
\left. \begin{cases} 
\; k_{\,11} \,, & M\,(\,t\,) \ < \ 1 \\
\; k_{\,12}\,, & M\,(\,t\,)\ \geqslant \ 1 
\end{cases}
\ \right. \,,\\
k_{\,2}\,(\,M\,) & \egal \max \Biggl\{\, 1 \moins \exp \biggl[\, 2.3 \cdot \biggl(\, M \,(\,t\,) \moins \Mmax\,(\,\phi\,) \, \biggl) \,\biggr] \,,\, 0\,\Biggr\} \,, 
\end{align*}
where the maximum mold growth value for these conditions is determined by: 
\begin{align}
\label{eq:Mmax}
\Mmax (\, \phi \,) & \egal A 
\plus B \, \biggl(\,\frac{\phic \moins \phi }{\phic \moins 100} \,\biggl)
\plus C \, \biggl(\,\frac{\phic \moins \phi}{\phic \moins 100} \,\biggl)^{\,2} \,,
\end{align}%
with $\phic$ (in percentage $\mathsf{\%}$) being the critical relative humidity to initiate mold growth. According to \cite{Viitanen2010}, four mold growth vulnerability classes of materials were defined with the corresponding values of the parameters $k_{\,1}\,$, $\phic\,$, $A\,$, $B\,$ and $C\,$ reported in Table~\ref{tb:mold_sensitivity_classes}.

\begin{table}
\centering
\begin{tabular}{l|cccccc}
\hline
\hline
\textit{Vulnerability classes} 
& $k_{\,11} $ 
&  $k_{\,12}$
& $A$
& $B$
& $C$ 
& $\phic \ \bigl[\,\mathsf{\%}\,\bigr]$ \\
\hline
\hline
\textrm{Very vulnerable}
& $1$
& $2$
& $1$
& $7$
& $-\,2$
& $80$ \\
\textrm{Vulnerable}
& $0.578$
& $0.386$
& $0.3$
& $6$
& $-\,1$
& $80$ \\
\textrm{Medium resistant}
& $0.072$
& $0.097$
& $0$
& $5$
& $-\,1.5$
& $85$ \\
\textrm{Resistant}
& $0.033$
& $0.014$
& $0$
& $3$
& $-\,1$
& $85$\\
\hline
\hline
\end{tabular}\bigskip
\caption{\small\em Mold growth vulnerability classes of materials \cite{Hukka1999, Viitanen2010}.}
\label{tb:mold_sensitivity_classes}
\end{table}


\subsection{Structural identifiability of the parameters}
\label{sec:struct_identifiability}

This section aims at justifying the identifiability of the unknown parameters: (i) $k_{11}$ the rate of increase when $M \ < \ 1$, (ii)  $k_{12}$ the rate of increase when $M \ \geqslant \ 1$ and (iii) the maximum mold growth index  $\Mmax$  defined through the coefficients $A\,$, $B\,$ and $C$. It corresponds to a total of five parameters:  
\begin{align*}
 \gP \eqdef \bigl\{\, P_{\,i}\,\bigr\} \egal \bigl\{\, k_{\,11} \,,\, k_{\,12} \,,\, A \,,\, B \,,\, C \,\bigr\} \,.
\end{align*}
It should be noted that parameter $k_{\,2}$ is not added since it is defined through parameters $\Mmax$ and therefore through $A\,$, $B\,$ and $C\,$.

A parameter $P_{\,i} \ \in \ \gP$ is Structurally Globally Identifiable (SGI) if the following condition is satisfied \cite{Walter1982}:
\begin{align*}
\forall t \,, \qquad M (\, P\,) \egal M (\, P^{\,\prime}\,) \quad \Rightarrow \quad P_{\,i} \egal P_{\,i}^{\,\prime} \,.
\end{align*}
It is assumed that $M$ is observable, that the function $f$ is known and that the derivative of $M$ is non-vanishing on any time interval, \ie $\od{M}{t} \ \neq \ 0$. Thus, the model writes: 
\begin{align*}
\od{M}{t} \egal \frac{k_{\,1}\,(\,t\,) \cdot k_{\,2} \,(\,t\,)}{ f} \,.
\end{align*}
So as to prove identifiability, it is assumed that another set of parameters, denoted with a superscript $^{\,\prime}\,$, holds:
\begin{align*}
\od{M^{\,\prime}}{t} \egal \frac{k_{\,1}^{\,\prime}\,(\,t\,) \cdot k_{\,2}^{\,\prime} \,(\,t\,)}{ f} \,.
\end{align*}
If $M\,(\,t\,) \ \equiv \ M^{\,\prime}\,(\,t\,)$, then $\od{M}{t} \ \equiv \ \od{M^{\,\prime}}{t} $ and therefore, since $\phi$ and $T$ can be measured, they are identical and we have:
\begin{align*}
k_{\,1}\,(\,t\,) \cdot k_{\,2} \ \equiv \ k_{\,1}^{\,\prime}\,(\,t\,) \cdot k_{\,2}^{\,\prime} \,.
\end{align*}
It results in:
\begin{align}
\label{eq:ident_k1}
k_{\,1}\,(\,t\,) & \cdot \max \Biggl\{\, 1 \moins \exp \biggl[\, 2.3 \, \biggl(\, M \,(\,t\,) \moins \Mmax\,(\,\phi\,) \, \biggl) \,\biggr] \,,\, 0\,\Biggr\}
\\ \ \equiv \
 & k_{\,1}^{\,\prime}\,(\,t\,) \cdot \max \Biggl\{\, 1 \moins \exp \biggl[\, 2.3 \, \biggl(\, M \,(\,t\,) \moins \Mmax^{\,\prime}\,(\,\phi\,) \, \biggl) \,\biggr] \,,\, 0\,\Biggr\} \nonumber \,.
\end{align}
It can be noted that $\Mmax$ corresponds to the asymptotic value of $M$:
\begin{align*}
\Mmax \egal \lim_{\,t \rightarrow \infty} M (\,t\,) \,.
\end{align*}
Thus, from Eq.~\eqref{eq:ident_k1} it can be deduced that:
\begin{align*}
\Mmax \ \equiv \ \Mmax^{\,\prime} \,,
\end{align*}
which yields $k_{\,1} \egal k_{\,1}^{\,\prime}$, ensuring that parameter $k_{\,1}$ is identifiable. According to its definition, at least one measurement for $M \ < \ 1$ and one for $M \ \geqslant \ 1$ should be enough to estimate this parameter. Moreover, from the equality $\Mmax \egal \Mmax^{\,\prime}$, we obtain:
\begin{align*}
A \plus B \, \biggl(\,\frac{\phic \moins \phi\,(\,t\,) }{\phic \moins 100} \,\biggl)
& \plus C \, \biggl(\,\frac{\phic \moins \phi\,(\,t\,) }{\phic \moins 100} \,\biggl)^{\,2} \\
& \ \equiv \ 
A^{\,\prime} 
\plus B^{\,\prime} \, \biggl(\,\frac{\phic \moins \phi\,(\,t\,) }{\phic \moins 100} \,\biggl)
\plus C^{\,\prime} \, \biggl(\,\frac{\phic \moins \phi\,(\,t\,) }{\phic \moins 100} \,\biggl)^{\,2} \,.
\end{align*}
As a consequence,
\begin{align*}
& A \egal A^{\,\prime} \,, 
&& B \egal B^{\,\prime} \,, 
&& C \egal C^{\,\prime} \,.
\end{align*} 
and the parameters are structurally globally identifiable. Therefore, it is necessary and sufficient to have measurements for at least three values of relative humidity $\phi$ to have the identifiability of parameters $A\,$, $B\,$ and $C\,$. One may conclude that if $M(\,t\,)$ is a solution of the model for a given set of parameters, and if $f$ is assumed to be known, then the parameters of the VTT model in $\gP$ are identifiable.


\subsection{Parameter estimation problem}

Since it was demonstrated that the unknown parameters $\gP$ are identifiable, it is important to detail the methodology to solve the parameter estimation problem, \ie to identify them. It is assumed that a set of Projected Observed Data (POD) of the mold growth index $M$ is available for a given relative humidity $\phi\,$:
\begin{align*}
\gM_{\mathrm{POD}} \bigl[ \, \phi \, \bigr] & \egal \bigl\{\, M_{\,n} (\, t_{\,i} \,) \,\bigr\} \,, 
&& n \in \bigl\{\, 1 \,,\, \ldots \,,\, N_{\,e} \,\bigr\} \,, 
 \qquad i \in \bigl\{\, 1 \,,\, \ldots \,,\, N_{\,t} \,\bigr\} \,,
\end{align*}
where $N_{\,e}$ is the number of experiments and $N_{\,t}$ the number of elements of the time grid.

The identification can be formulated as an optimization problem. The estimated parameter, denoted with the super script $^{\,\circ}\,$, is determined according to:
\begin{align*}
  \gP^{\,\circ} \egal \arg \min_{\,\gP} \bigl(\, \mathrm{J} (\, \gP \,) \,\bigr) \,,
\end{align*}
where the cost function $\mathrm{J}$ is defined as:
\begin{align}\label{eq:cost_function}
  \mathrm{J} (\, \gP \,) \egal \Bigl|\Bigl|\, \gM_{\,\mathrm{POD}} \moins \mathcal{P} \bigl(\, M (\, \gP\,)\,\bigr) \,\Bigr|\Bigr|^{\,2} \,,
\end{align}
where $\mathcal{P}$ is used to project the numerical solution $M$ from Eq.~\eqref{eq:odeM} on the time grid of the discrete observations $\gM_{\,\mathrm{POD}}$. For discrete quantities, the $\mathscr{L}_{\,2}$ norm is usually defined as: 
\begin{align*}
  \Bigl|\Bigl|\, u \moins v \,\Bigr|\Bigr|^{\,2} \egal \sum_{i \, = \, 1}^{N_{\,t}} \bigr(\, u\,(\,t_{\,i}\,) - v\,(\,t_{\,i}\,) \bigl)^{\,2} \,,
\end{align*}
where $N_{\,t}$ is the number of temporal time steps.

The residual between the index $M$ computed with the numerical model and the experimental data is defined as: 
\begin{align*}
\varepsilon\,(\,t\,) \eqdef \Bigl|\Bigl|\, \gM_{\,\mathrm{POD}} \moins \mathcal{P} \bigl(\, M (\, \gP^{\,\circ}\,)\,\bigr) \,\Bigr|\Bigr| \,.
\end{align*}
To estimate the quality of the solution of the parameter identification, the normalized \textsc{Fisher} matrix \cite{Karalashvili2015, Ucinski2004} is defined according to: 
\begin{subequations}
\label{eq:fisher_matrix}
\begin{align}
 \gF & \egal \bigl[\, F_{\,i \, j} \,\bigr] \,, && \forall \ (\,i\,,\,j\,) \ \in \ \bigl\{ 1, \ldots, N_{\,P} \, \bigr\}^{\,2}\,, \\[3pt]
 F_{\,i \, j} & \egal \sum_{\,n \,= \,1}^{\,N_{\,e}} \ \int_{\,0}^{\, \tau} \Theta_{\,i \,n} \ (t) \ \Theta_{\,j \,n} \ (t) \ \mathrm{dt} \,,
\end{align}
\end{subequations}
where  $N_{\,p}$ is the number of parameters and  $\Theta_{\,i \,n}$ is the sensitivity coefficient of the solution $M$ upon the parameter $P_{\,i}$ computed for the experiment number $n\,$. The sensitivity coefficient is defined as \cite{Finsterle2015}: 
\begin{align*}
& \Theta_{\,\,i \,n} \,(t) \egal
\frac{\sigma_{\,p}}{\sigma_{\,M}} \, \pd{M_{\,n}}{P_{\,i}} \,, 
&& \forall \ i \ \in \   \bigl\{ 1, \ldots, N_{\,p} \, \bigr\} \,,
\end{align*}
where $\sigma_{\,M}$ and $\sigma_{\,p}$ are scaling factors in order to define dimensionless sensitivity coefficients. Since we are using dimensionless variables, these factors are set to unity $\sigma_{\,M} \egal \sigma_{\,p} \egal 1 \,$.

The function $\Theta$ is computed by taking the derivative of Eq.~\eqref{eq:odeM} with respect to the required parameter. The matrix $\gF$ measures the total sensitivity of the system for the measurements \cite{Nenarokomov2005,Artyukhin1985}, to variations of the entire set of parameters $\gP\,$. Under some assumptions (see \cite{Walter1990}), its inverse is the matrix of variance of the parameters considered as random variables of the observable fields. It summarizes the quality of the information obtained in the parameter identification process \cite{Guisasola2006}. It can be used to assess the estimation uncertainty by computing an error estimator for the parameter $P_{\,i}\,$:
\begin{align*}
  \eta_{\,i} \egal \sqrt{ (\,\gF^{\,-1}\,)_{\,i\,i}} \,,
\end{align*}
A high value of $\eta_{\,i}$ indicates a potentially high error in the parameter identification of $P_{\,i} \,$. It should be noted that the matrix $\mathds{F}^{\,-1}$ provides the lower bound for the asymptotic co-variance matrix, assuming no data auto-correlation, white noise measurement or uncorrelated errors \cite{Walter1990}.


\section{Experimental data}
\label{sec:experimental_data}

This Section presents the experimental design to obtain observations of mold growth. Once collected, it will be shown that the discrete observation data need to be projected in order to proceed to parameter identification. 


\subsection{Experimental facility}

\begin{figure}
\begin{center}
\subfigure[]{\includegraphics[width=.55\textwidth]{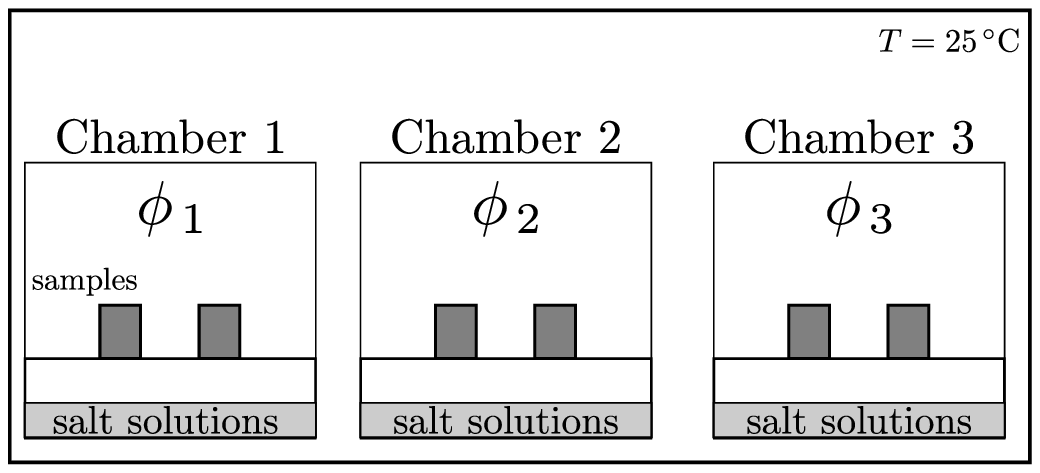}} \hspace{0.2cm}
\subfigure[]{\includegraphics[width=.25\textwidth]{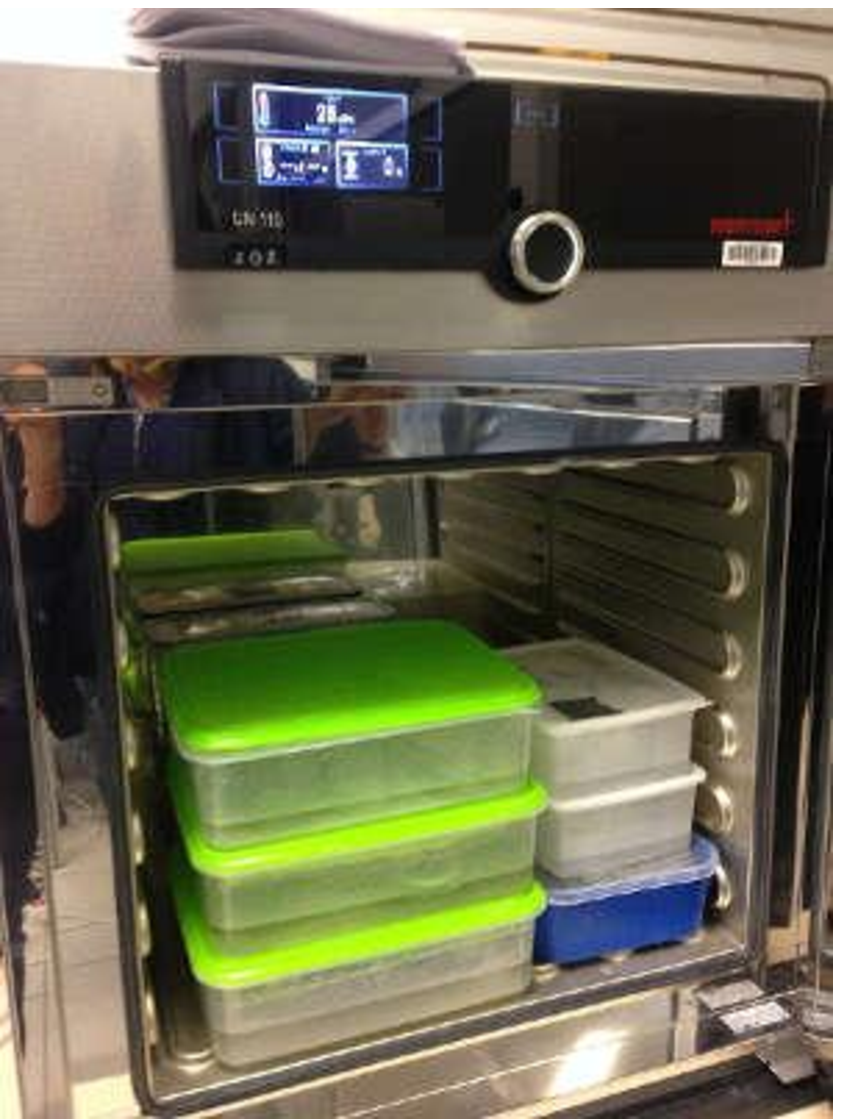}}
\caption{\small\em Illustration (a) and picture (b) of the experimental facility.}
\label{fig:facility}
\end{center}
\end{figure}

The material under investigation is bamboo fiberboard manufactured without glue. It can be used for wall insulation in buildings. The bamboo belongs to the \emph{Bambusa stenostachya} species and interested readers may refer to \cite{Nguyen2017} for more information on this material. To investigate the critical moisture levels of the fiberboard, the test specimens are cut into samples measuring $50 \times 50 \times 6 \ \mathsf{mm}\,$. The test facility used to provide the experimental data is composed of three climatic chambers as illustrated in Figure~\ref{fig:facility}. Two samples of the same species are settled in each chamber. For each chamber, the relative humidity is controlled using saturated salt solutions. The samples are settled above the solution on grids. Three experiments are conducted, corresponding to three different relative humidity values:
\begin{align*}
  & \phi_{\,1} \egal 0.75 \,, && \phi_{\,2} \egal 0.84 \,, && \phi_{\,3} \egal 0.97 \,.
\end{align*}

The three chambers are located in an oven so as to control the temperature, set at $T \egal 25 \ \mathsf{^\circ C}$ for the study. Wireless sensors inside the small chambers are used to check the temperature and relative humidity levels. There were no spore inoculation during the experimental process. Thus, the study rely only on the spores present in the air.

The VTT model takes into account the surface material temperature and relative humidity. Measurements are made in the chamber, assuming the measures are the same as at the material's surface. If we assume also that the heat and moisture diffusion characteristic times are smaller than the measurement time step. This assumption is verified by evaluating the characteristic time of heat and moisture diffusion and the \textsc{Fourier} number for moisture diffusion given by:
\begin{align*}
  & \mathrm{Fo} \egal \frac{\delta_{\,v} \; t^{\,\circ}}{\xi \; L^{\,2}} \,,
\end{align*}
where $\delta_{\,v}$ is the vapor permeability, $\xi$ is the sorption curve, $L$ the characteristic length and $ t^{\,\circ}$ the characteristic time of the problem. This number reflects the importance of the moisture diffusion in the material. It can also be seen as a ratio between the characteristic time of the problem and of the moisture diffusion:
\begin{align*}
& \mathrm{Fo} \egal \frac{t^{\,\circ}}{t^{\,\mathrm{d}}} \,,
\end{align*}
where the characteristic moisture diffusion time $t^{\,\mathrm{d}}$ equals:
\begin{align*}
t^{\,\mathrm{d}} \egal \frac{\xi \; L^{\,2}}{\delta_{\,v}} \,.
\end{align*}
According to \cite{Nguyen2017}, the vapor permeability and the sorption capacity of the bamboo fiber at $T \egal 25 \ \mathsf{^{\,\circ}C}\,$ scale with: 
\begin{align*}
  & \delta_{\,v} \ \approx \ 2.7 \cdot 10^{\,-16} \ \bigl[\, \mathsf{h} \,\bigr] \,,
&& \xi \egal \frac{1}{P_{\,\mathrm{sat}}\,(\,T\,)} \; \pd{w}{\phi} \ \approx \ 4 \cdot 10^{\,-5} \ \bigl[\, \mathsf{kg/(m^{\,3}.Pa)} \,\bigr] \,.
\end{align*}
Knowing that $L \egal 0.05 \ \mathsf{m} \,$, the characteristic time of moisture diffusion is approximated by:
\begin{align*}
  t^{\,\mathrm{d}} \ \approx \ 3 \ \bigl[\,\mathsf{h}\,\bigr] \,.
\end{align*}
This means that if a solicitation of relative humidity occurs at a boundary, the moisture process will take around $3 \ \mathsf{h}$ to reach the other boundary. Since the experimental assessment takes several weeks, the assumption is verified. Since the temperature diffusion in the material is faster than the moisture diffusion, the assumption is also verified for the temperature measurements.

\subsection{Mold growth data}
\label{sec:mold_growth_data}

\begin{figure}
\begin{center}
\subfigure[$t \egal 0 \ \mathsf{days} \,,\, M \egal 0$]{\includegraphics[width=.30\textwidth]{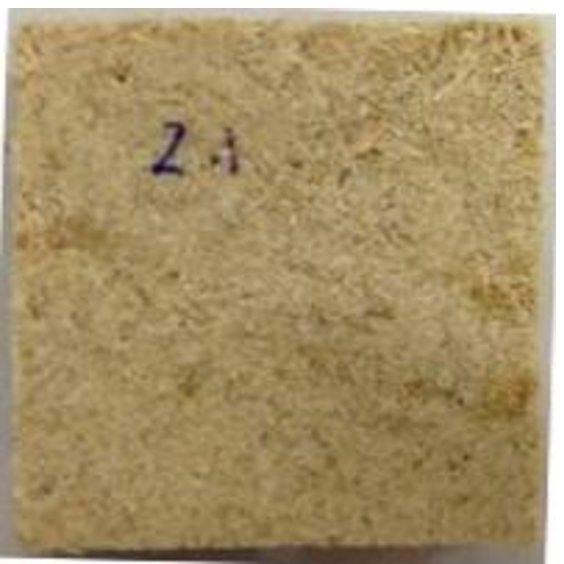}} \hspace{0.2cm}
\subfigure[$t \egal 8 \ \mathsf{days} \,,\, M \egal 4$]{\includegraphics[width=.30\textwidth]{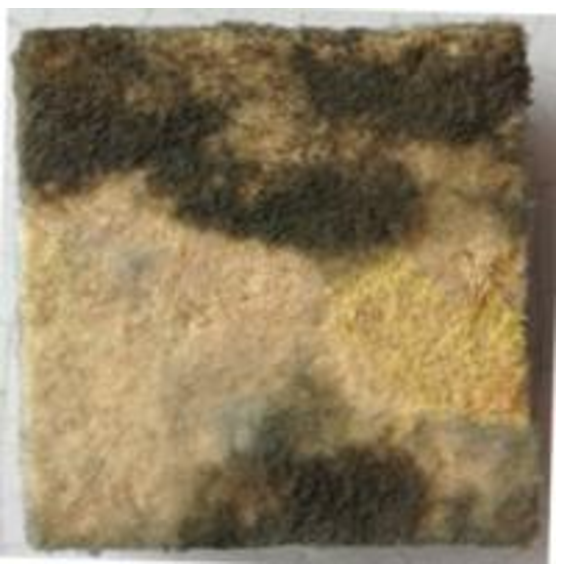}}
\hspace{0.2cm}
\subfigure[$t \egal 15 \ \mathsf{days} \,,\, M \egal 6$]{\includegraphics[width=.30\textwidth]{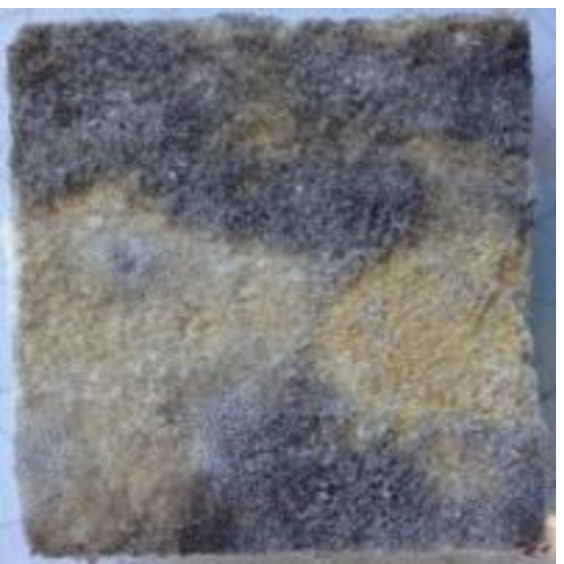}}
\caption{\small\em Illustration of mold growth on the fiberboard for $\phi \egal 0.97 \,$ at different times.}
\label{fig:samples_pictures}
\end{center}
\end{figure}

The investigation was conduced over $16$ weeks with daily assessment of mold growth. Mold was observed using an optical microscope (Leica DM4000M) with $5$--$40$ magnification and the naked eye. The experimental data were evaluated according to the rating scales of mold assessment defined in \cite{Hukka1999,Viitanen2010} and recalled in Table~\ref{tb:Mold_index}. Figure~\ref{fig:samples_pictures} shows the mold growth on the fiberboard for $\phi_{\,3} \egal 0.97 \,$.

\begin{table}
\centering
\begin{tabular}{c|l}
\hline
\hline
\textit{Index} & \textit{Description of the growth rate}\\
\hline
\hline
0 & No mold growth \\
1 & Small amounts of mold growth detected with microscopy \\
2 & Moderate mold growth detected with microscopy (coverage more than $10 \%$) \\
3 & Some growth detected visually \\
4 & Visually detected coverage more than $10 \%$\\
5 & Visually detected coverage more than $50 \%$\\
6 & Coverage about $100 \%$\\
\hline
\hline
\end{tabular}\bigskip
\caption{\small\em Mold index according to \cite{Hukka1999, Viitanen2010}.}
\label{tb:Mold_index}
\end{table}

The discrete experimental data corresponding to three experiments are shown in \mbox{Figure~\ref{fig:Dataexp_proj1}}. The data are discontinuous and the index $M$ is only evaluated as an integer. It should be noted that both samples give exactly the same experimental observations. For $\phi_{\,1} \egal 0.75 \,$, mold required almost $90\ \mathsf{days}$ to grow. It appears faster for $\phi_{\,2} \egal 0.84\,$, around $15 \ \mathsf{days}$ and for $\phi_{\,3} \egal 0.97 \,$, around $7 \ \mathsf{days}\,$. Considering these observations, it is reasonable to estimate that the bamboo fiberboard belongs to the very vulnerable class of materials. 


\subsection{Projected mold growth data}
\label{sec:projected_mold_growth_data}

In the first step, we show that the discrete measured quantities may not fit the continuous model defined by Eq.~\eqref{eq:odeM}. Subsequently, this analysis demonstrates the need to project the discrete observed quantities on defined continuous functions. In the second step, the parameter estimation problem is solved using the POD.


\subsubsection{Need to project the discrete experimental data}

According to the definition of the observable index $M$ recalled in Table~\ref{tb:Mold_index}, it is a difficult task to define intermediate values in order to provide more experimental data. Moreover, looking at the experimental measurements, illustrated for the second experiment in Figure~\ref{fig:exp2_proj1}, for a constant $\phi\,$, the discrete observed data of $M$ remains vanishing for $t \ \in \ \bigl[\,0 \,,\, 13 \,\bigr] \ \mathsf{days}\,$. Then, assuming the continuous model and the function $f$ are determined, $\od{M}{t} \egal 0$ yields:
\begin{align*}
  1 \moins \exp \Bigl[\, 2.3 \cdot \bigl(\, M \,(\,t\,) \moins \Mmax\,(\,\phi\,) \, \bigl) \,\Bigr] \ \leqslant \ 0 \,,
\end{align*}
which is equivalent to:
\begin{align}\label{eq:MinfMmax}
  M\,(\,t\,) \ \geqslant \ \Mmax \,.
\end{align}
Since we assume $M$ to vanish and $M$ to tend in an increasing way to $\Mmax \ > \ 0\,$, Eq.~\eqref{eq:MinfMmax} is impossible. There is no time $t$ verifying this equation. Therefore, the continuous model cannot match the experimental data.


\subsubsection{Projection of the discrete experimental data}
\label{sec:projection_ODD_VTT}

\begin{figure}
\centering
\subfigure[$\phi_{\,1} \egal 0.75$ \label{fig:exp1_proj1}]{\includegraphics[width=0.45\textwidth]{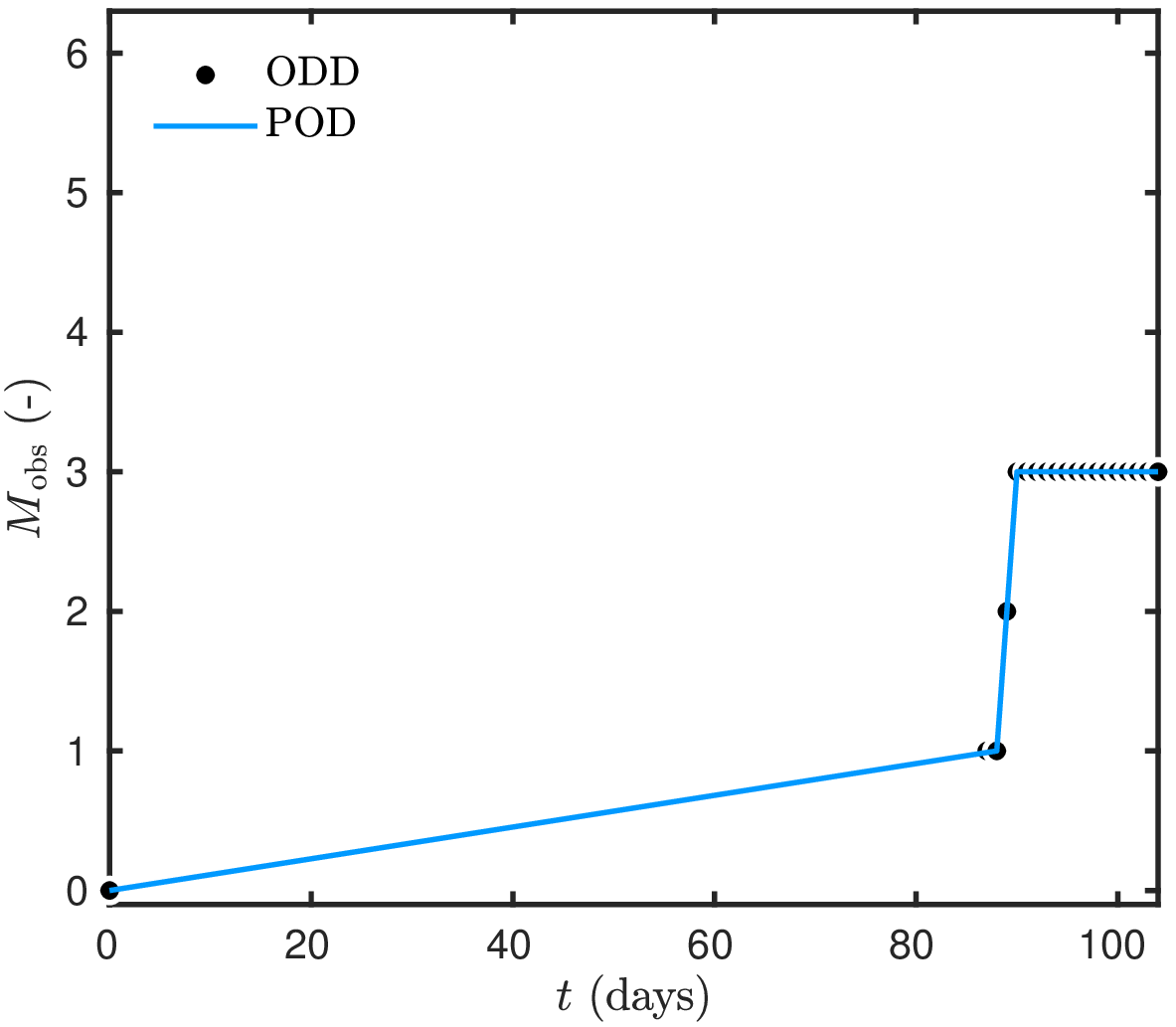}} \hspace{0.2cm}
\subfigure[$\phi_{\,2} \egal 0.84$ \label{fig:exp2_proj1}]{\includegraphics[width=0.45\textwidth]{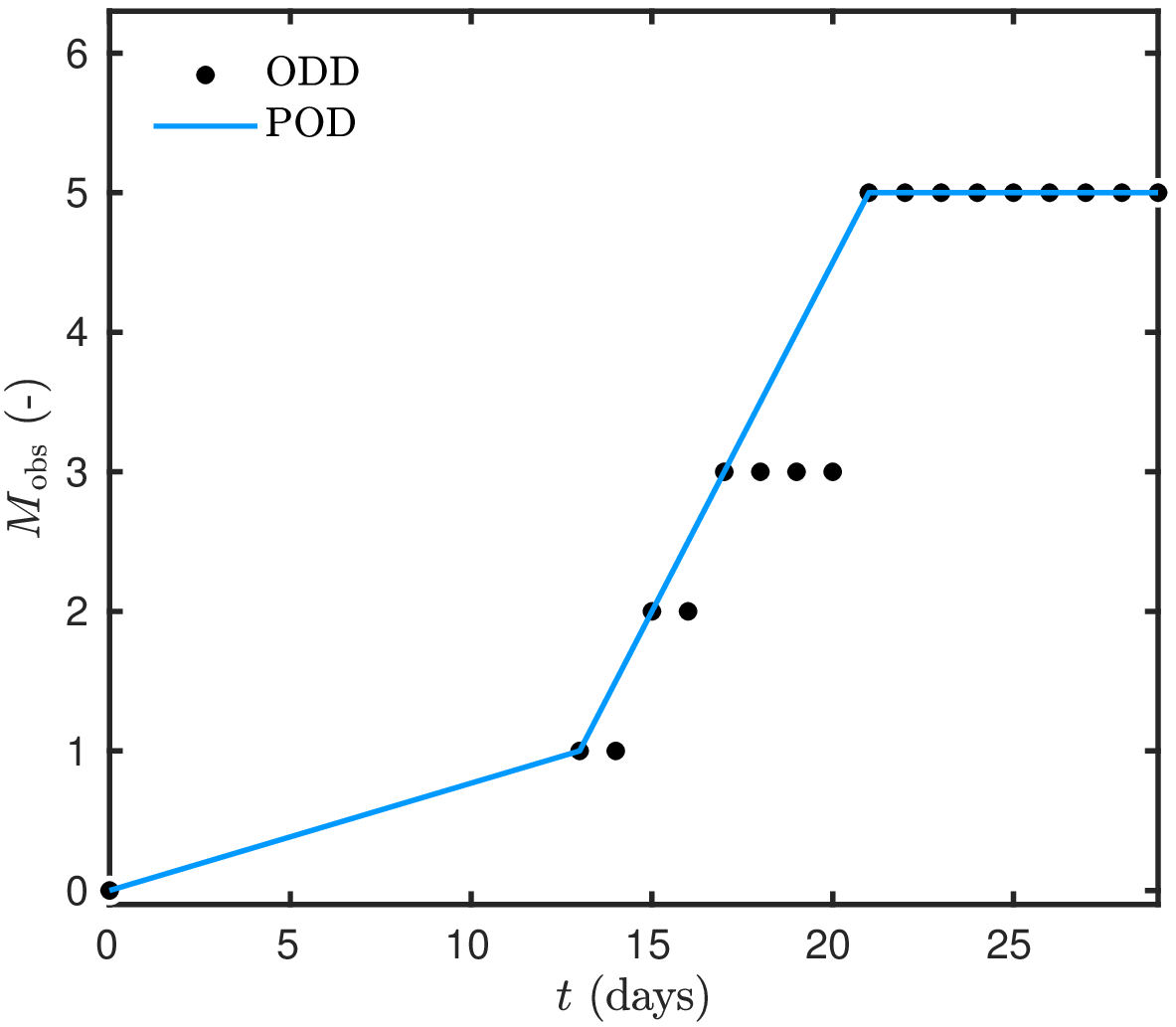}} \\
\subfigure[$\phi_{\,3} \egal 0.97$ \label{fig:exp3_proj1}]{\includegraphics[width=0.45\textwidth]{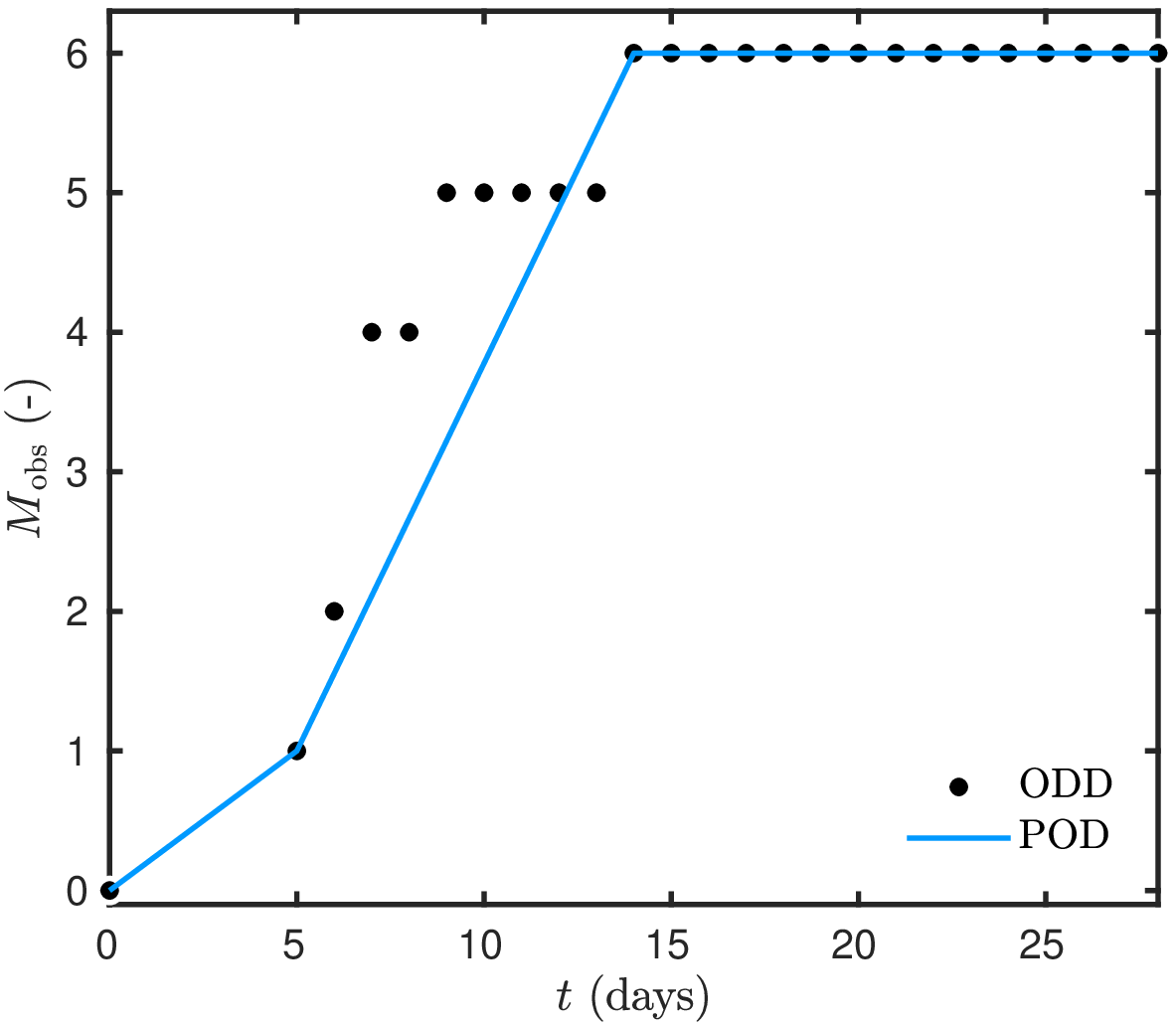}}
\caption{Experimental ODD available for three different values of relative humidity and the POD resulting from the projection.}
\label{fig:Dataexp_proj1}
\end{figure}

To solve this problem, the experimental data are projected on the following intermediate family of functions:
\begin{align*}
  g_{\,a,\,b,\,c\,}(\,t\,) \egal & 
  \left. \begin{cases} 
  \; \dfrac{1}{a} \,t \,, & t \ < \ a \,,\\[10pt]
  \; \dfrac{b \moins 1}{c \moins a} \, \bigl(\,t \moins c \,\bigr) \plus b \,, & t \ \in \ \bigl[\, a \,,\, c \,\bigr] \,,\\[10pt]
  \; b \,, & t \ > \ c \,,
\end{cases}
\ \right. \, \quad \text{with } a \ < \ b \ < \ c \,.
\end{align*}
Function $g$ is parametrized using three parameters: (i) the time $a$ to reach $M \egal 1\,$, (ii) the maximum value of $M\,$, denoted $b\,$, and (iii) the time $c$ to reach this maximum value:
\begin{align*}
  & a \egal \arg_{\,t} \bigl(\, M(\,t\,) \egal 1 \,\bigr) \,, 
 && b \egal \max \bigl(\, M(\,t\,) \,\bigr) \,,
 && c \egal \arg_{\,t} \bigl(\, M(\,t\,) \egal b \,\bigr) \,. 
\end{align*}
The values of parameters $(\,a\,,\,b\,,\,c\,)$ are given in Table~\ref{tb:coeff_proj}. The time domain is discretized with a uniform mesh using a time step of $6 \ \mathsf{min}\,$. Figure~\ref{fig:Dataexp_proj1}(a,b,c) shows the comparison between the discrete and projected experimental data, which is satisfactory. We can use these projected data to identify the parameters of the VTT model.

\begin{table}
\centering
\begin{tabular}{c|c|c|c}
\hline
\hline
\textit{Parameters} & \textit{Experiment} $\phi_{\,1}$ & \textit{Experiment} $\phi_{\,2}$ & \textit{Experiment} $\phi_{\,3}$\\
\hline
\hline
$a$ & $88$ & $13$ & $5$ \\
$b$ & $3$ & $5$ & $6$ \\
$c$ & $90$ & $21$ & $14$ \\
\hline
\hline
\textit{$\L_{\,2}$ residual}
& $0.002$ & $0.47$ & $0.65$   \\
\hline
\hline
\end{tabular}\bigskip
\caption{\small\em Parameters for the projection of the discrete experimental data using the family of functions $g_{\,a,\,b,\,c\,}(\,t\,)$ to obtain the projected observed data.}
\label{tb:coeff_proj}
\end{table}

\section{Estimating the VTT model parameters and their sensitivity}
\label{sec:estimation_parameters_VTT}

Starting form our POD, we have a twofold issue. Our first goal is to identify first $\bigl(\, k_{\,11} \,,\, k_{\,12} \,,\, \Mmax \,\bigl)$ and subsequently $\bigl(\, A \,,\, B \,,\, C\,\bigr) \,$. The second goal is to estimate the quality of this identification. This last results enables to quantify the robustness of the VTT model.

For the sake of clarity, it is important to distinguish: (i) the Observed Discrete Data (ODD) corresponding to the experimental data, (ii) the Projected Observed Data (POD) resulting from the projection of the experimental data on defined functions and (iii) the Simulated Quantities (SQ) computed from the mathematical model of mold growth. Different types of simulated quantities can be distinguished. On the one hand, the simulated quantities can be obtained with parameters found in the literature and reported in Table~\ref{tb:mold_sensitivity_classes}. Analyzing the experimental results, and particularly Figure~\ref{fig:samples_pictures}, the material seems very vulnerable to mold growth. As a consequence, these parameters will be those of the corresponding class in Table~\ref{tb:mold_sensitivity_classes}. They are used as an initial guess in the optimization procedure. These parameters are denoted as \emph{a priori}, in the sense that they are suggested by the VTT model. Using the values of parameters $A\,$, $B\,$ and $C\,$, the \apriori value of parameter $\Mmax$ can be computed using Eq.~\eqref{eq:Mmax}. The critical relative humidity is set to $\phi_{\,c} \egal 0.75$ since mold growth was observed for this condition. On the other hand, the simulated quantities can be computed with the estimated parameters, which are the solution of the optimization process.


\subsection{Estimation of parameters $\bigl(\, k_{\,11} \,,\, k_{\,12} \,,\, \Mmax \, \bigr)$}

The first step consists in defining three cost functions according to three experiments:
\begin{align*}
  & \mathrm{J}^{\,n} (\, \gP \,) \egal \Bigl|\Bigl|\, \gM_{\,\mathrm{POD}}\,(\,\phi_{\,n}\,) \moins \mathcal{P} \bigl(\, M (\, \gP\,)\,\bigr) \,\Bigr|\Bigr|^{\,2} \,, 
  && n \in \bigl\{\, 1 \,,\, \ldots \,,\, 3 \,\bigr\} \,. 
\end{align*}

Each cost function is minimized using the least squares method. The identification of the set of parameters $\bigl(\, k_{\,11} \,,\, k_{\,12} \,,\, \Mmax \, \bigr)$ is realized using  the \texttt{fmincon} from the \texttt{Optimization Toolbox} in the \texttt{Matlab\texttrademark}. It provides an efficient interiorpoint algorithm with (linear or nonlinear) constraints on the unknown parameters \cite{Byrd2000}. Considering the physical problem, the only obvious constraint is on parameter $\Mmax$ by the definition of the mold growth index in Table~\ref{tb:Mold_index}:
\begin{align*}
  \Mmax \ \in \ \bigl[ \, 0 \,, 6 \, \bigl] \,.
\end{align*}

As a result of the optimization process, the estimated values are reported in Table~\ref{tb:unknown_parameters}. The POD and the SQ match well while the \apriori parameters were poor.


\subsection{Estimation of parameters $\bigl(\, A \,,\, B \,,\, C \, \bigr)$}

As mentioned in Section~\ref{sec:struct_identifiability}, since measurements are available for three relative humidity values, it is possible to identify the parameters $A\,$, $B\,$ and $C\,$. These parameters are also estimated using the least square method. Here, we have a linear system of equations that can be written as:
\begin{align}\label{eq:system_eq_ABC}
  Y \egal S \, X \,,
\end{align}
with $Y$ being the vector containing the previously estimated parameters $\Mmax$ for the three experiments:
\begin{align*}
  Y \egal \bigl[\, \Mmax(\,\phi_{\,1}\,) \,,\, \Mmax(\,\phi_{\,2}\,) \,,\, \Mmax(\,\phi_{\,3}\,) \,\bigr]^{\,T} \ ,
\end{align*}
while the vector $X$ contains the unknown parameters: 
\begin{align*}
  X \egal \bigl[\, A \,,\, B \,,\, C \,\bigr]^{\,T} \,.
\end{align*}
The matrix $S$ is built using Eq.~\eqref{eq:Mmax}\,:
\begin{align*}
  S \egal
  \begin{bmatrix}
  1 
  & \biggl(\,\frac{\phic \moins \phi_{\,1} }{\phic \moins 1} \,\biggl) 
  & \biggl(\,\frac{\phic \moins \phi_{\,1}}{\phic \moins 1} \,\biggl)^{\,2} \\
  1 
  & \biggl(\,\frac{\phic \moins \phi_{\,2} }{\phic \moins 1} \,\biggl) 
  & \biggl(\,\frac{\phic \moins \phi_{\,2}}{\phic \moins 1} \,\biggl)^{\,2} \\
  1 
  & \biggl(\,\frac{\phic \moins \phi_{\,3} }{\phic \moins 1} \,\biggl) 
  & \biggl(\,\frac{\phic \moins \phi_{\,3}}{\phic \moins 1} \,\biggl)^{\,2}
  \end{bmatrix} \,.
\end{align*}

Thus, the system from Eq.~\eqref{eq:system_eq_ABC} is solved directly to determine parameters $\bigl(\, A \,,\, B \,,\, C \, \bigr)\,$. The estimated parameters are reported in Table~\ref{tb:unknown_parameters}. For comparison, the values corresponding to a very vulnerable materials are also provided. Since only three measurements are available, it is not possible to compute the confidence interval for the estimated parameters.

\begin{table}
\centering
\begin{tabular}{c|c|ccc}
\hline
\hline
\textit{Experiment} & \textit{Parameters} & $k_{\,11}$ & $k_{\,12}$ & $\Mmax$ \\[3pt]
\hline
\hline
\multirow{3}{*}{$\phi_{\,1} \egal 0.75$} & 
\emph{A priori} values & 1 & 2 & 1  \\
& Estimated values & 4.25 & 199 & 3.05  \\
& Error estimator $\eta$
& $1.54$
& $0.5$
& $0.018$ \\[3pt]
\hline
\multirow{3}{*}{$\phi_{\,2} \egal 0.84$} & \emph{A priori} values  & 1 & 2 & 3.26  \\
& Estimated values & 5.89 & 39.3 & 5.09   \\
& Error estimator $\eta$
& $3.92$
& $0.47$
& $0.096$ \\[3pt]
\hline
\multirow{3}{*}{$\phi_{\,3} \egal 0.97$} & \emph{A priori} values  & 1 & 2 & 5.61   \\
& Estimated values & 2.19 & 6.99 & 5.99   \\
& Error estimator $\eta$
& $6.00$
& $0.38$
& $0.071$ \\[3pt]
\hline
\multirow{2}{*}{All three}
& \emph{A priori} values & $1$ & $2$ & - \\
& Estimated values & $4.56$ & $198$ & - \\
& Error estimator $\eta$
& $3.05$
& $1.9$
& - \\[3pt]
\hline
\hline
\textit{Experiment} & \textit{Parameters} & $A$ & $B$ & $C$ \\[3pt]
\hline
\hline
\multirow{2}{*}{All three}
& \emph{A priori} values & $1$ & $7$ & $-2$  \\
& Estimated values & $3.05$ & $7.39$ & $-4.5$ \\[3pt]
\hline
\hline
\end{tabular}\bigskip
\caption{\small\em Estimated parameters of the VTT model for the bamboo woodfibers.}
\label{tb:unknown_parameters}
\end{table}


\subsection{Discussion of the results}

\begin{figure}
\centering
\subfigure[$\phi_{\,1} \egal 0.75$ \label{fig:M1_ft}]{\includegraphics[width=0.45\textwidth]{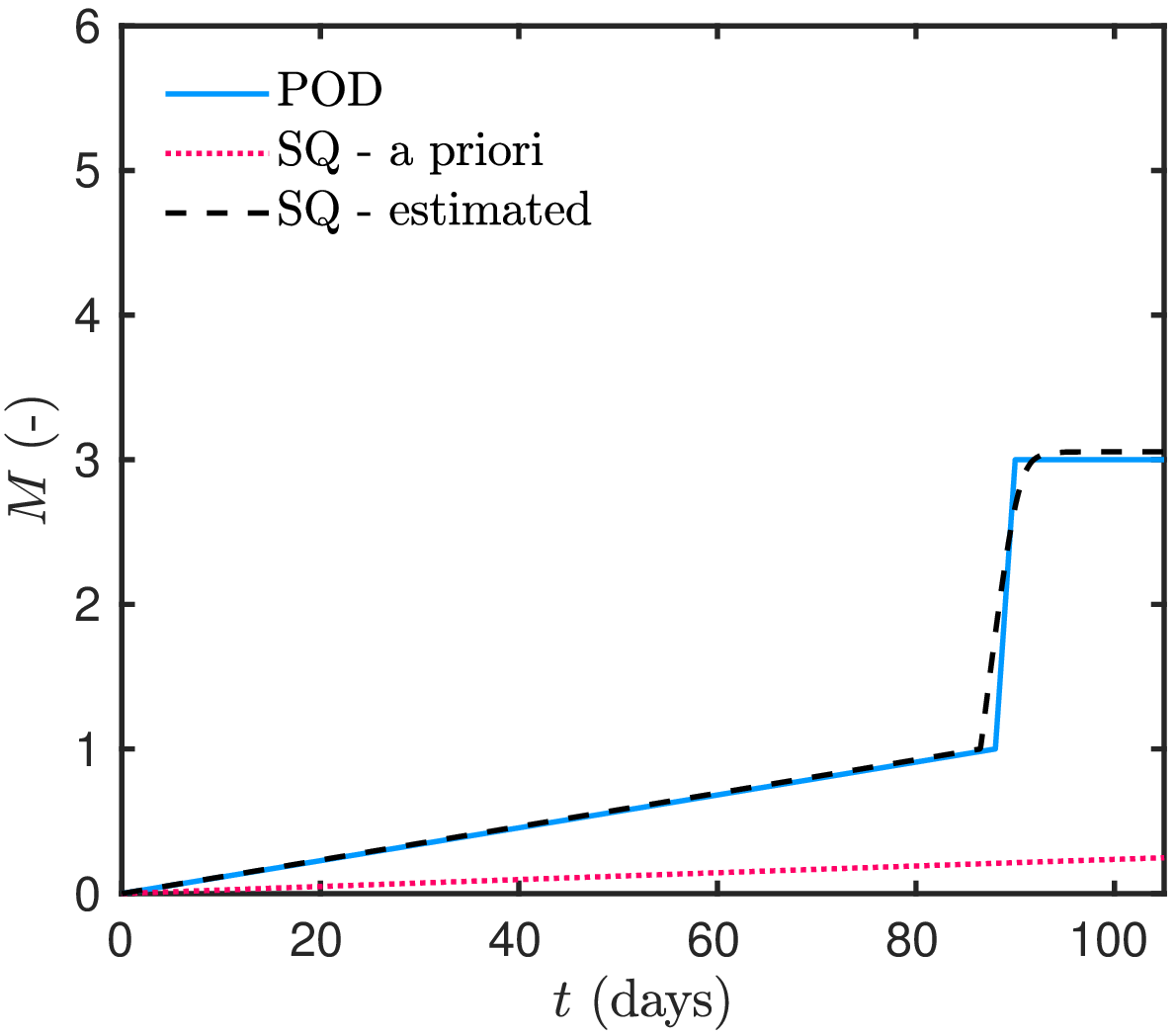}} \hspace{0.2cm}
\subfigure[$\phi_{\,2} \egal 0.84$ \label{fig:M2_ft}]{\includegraphics[width=0.45\textwidth]{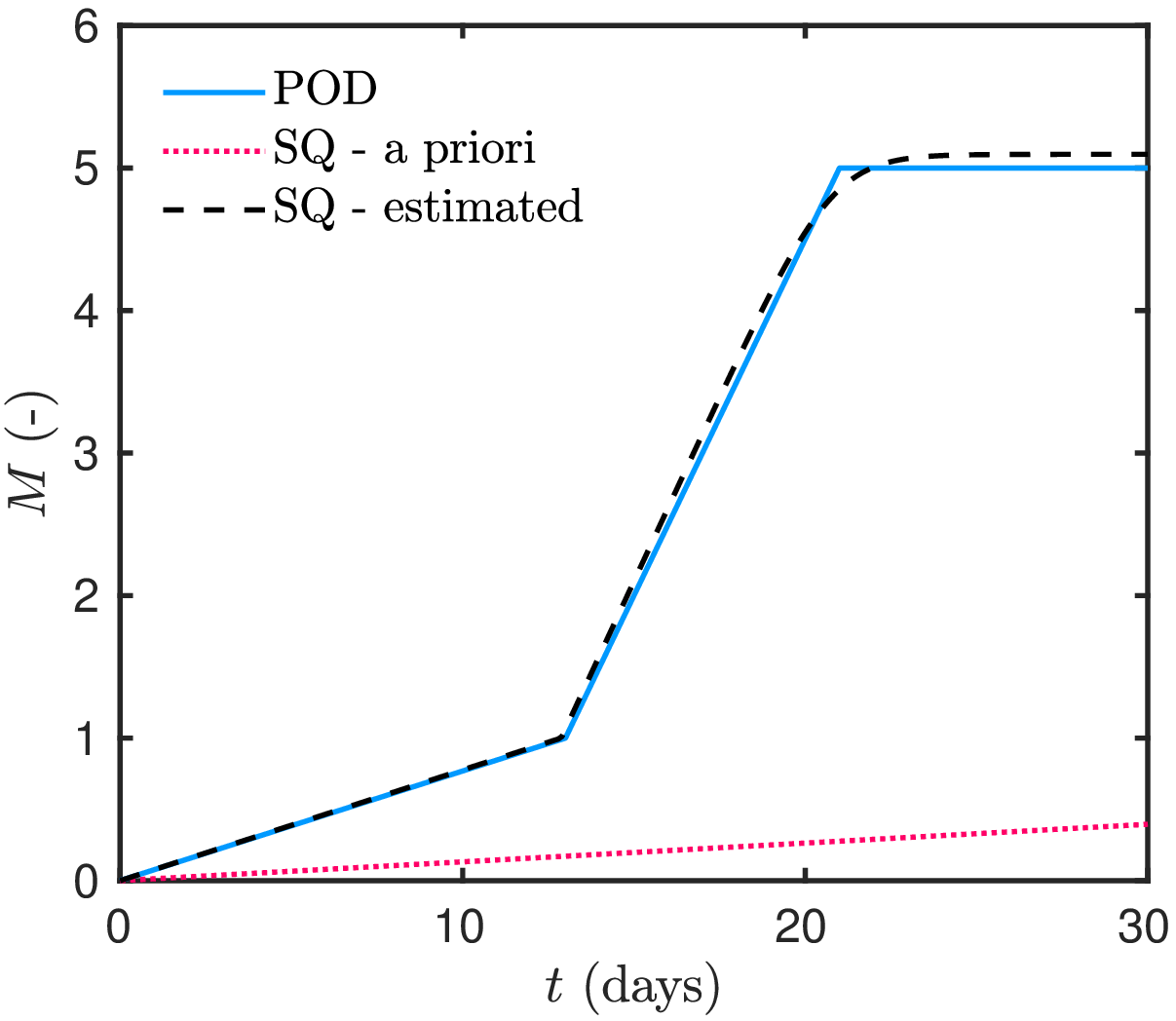}} \\
\subfigure[$\phi_{\,3} \egal 0.97$ \label{fig:M3_ft}]{\includegraphics[width=0.45\textwidth]{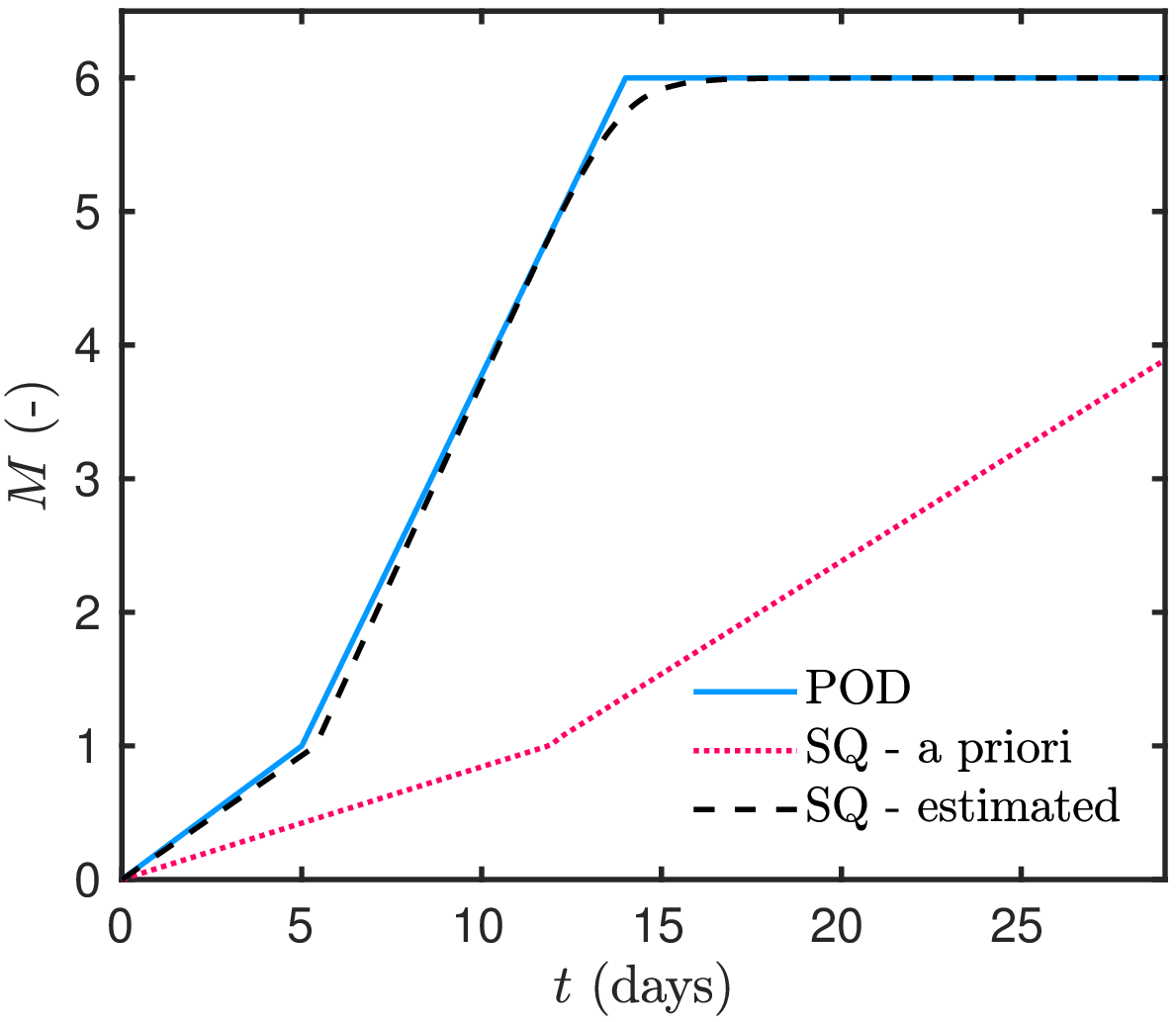}}\hspace{0.2cm}
\subfigure[\label{fig:res_ft}]{\includegraphics[width=0.45\textwidth]{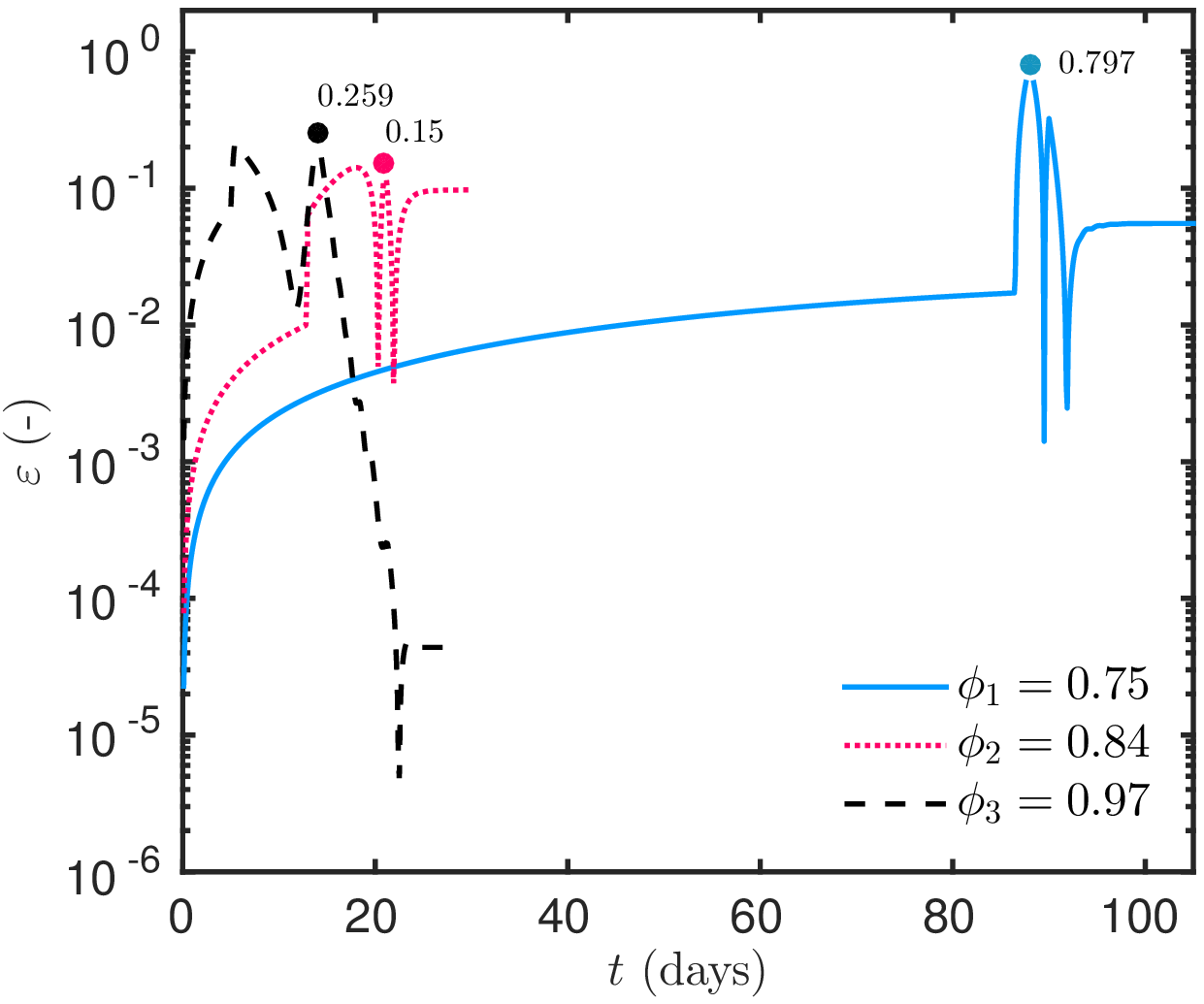}} 
\caption{\small\em \emph{(a,b,c)} Comparison of the POD with the SQ using the VTT model with the \apriori and estimated parameters. \emph{(d)} Residual between the POD and the SQ with the estimated parameters.}
\label{fig:M_ft}
\end{figure}

With the estimated parameters, a good agreement is observed between the numerical results and the projected observed data as shown in Figure~\ref{fig:M_ft}(a,b,c). The residual between measurements and the results of the numerical model is shown in Figure~\ref{fig:res_ft}. It remains at the order $\O(\,10^{\,-2}\,) \,$, which appears satisfactory. However, the residual increases at a particular moment, such as $t \egal 90 \ \mathsf{days}$, corresponding to sharp changes in the POD. The residual is not very satisfactory for these times. In addition, the $\L_{\,2}$ of the residual is $0.09\,$, $0.05$ and $0.05$ for each experiment.

As we can notice in Table~\ref{tb:unknown_parameters}, the estimated parameters are higher than \apriori values provided by the original model, particularly for the experiments with $\phi_{\,1}$ and $\phi_{\,2}\,$. As shown in Figure~\ref{fig:M_ft}(a,b,c), the mold index computed with \apriori parameter values underestimates the experimental data. In other words, the parameters of the most vulnerable material class are not able to represent the physical phenomena observed in these experiments. Moreover, according to the vulnerability classes recalled in Table~\ref{tb:mold_sensitivity_classes}, the parameters vary within intervals of variations $k_{\,11} \ \in \ \bigl[\, 0.033 \,,\, 1 \,\bigr]$ and $k_{\,12} \ \in \ \bigl[\, 0.014 \,,\, 2 \,\bigr]\,$. In our study, the estimated parameters $k_{\,11}$ and $k_{\,12}$ are completely out of this interval. Even if the discrepancies between the numerical results and the experimental data are low, these intermediate results highlight that the estimated parameters $\bigl(\,k_{\,11} \,,\, k_{\,12}\,\bigr)$ do not agree with classification of the original definition of the VTT model (Table~\ref{tb:mold_sensitivity_classes}).

Moreover, it is of major importance to quantify the confidence one may have in the estimated parameters. For this, the \textsc{Fisher} matrix is computed for parameters $k_{\,11}\,$, $k_{\,12}$ and $\Mmax $ according to Eq.~\eqref{eq:fisher_matrix} and the error estimators are reported in Table~\ref{tb:unknown_parameters}. The quality of the estimation is not satisfactory since the estimator $\eta$ is of the order $\O\,(1)\,$ for parameter $k_{\,11}$. The uncertainty on this parameter is very high. It corresponds to an uncertainty of $300\%$ on this parameter for the experiments $\phi_{\,3}\,$.

Further investigations are carried out to see if the estimation of $k_{\,11}$ and $k_{\,12}$ can be improved. Since these two parameters do not depend on the relative humidity nor on the temperature, they have been estimated with the three experiments together by defining one single cost function. The results are reported in Table~\ref{tb:unknown_parameters}. The $\L_{\,2}$ of the residual is $0.45\,$, $0.53$ and $2.15$ for each experiment. Thus, the quality of the estimation is not satisfactory and even worse than for the estimation using each experiment separately.

These intermediate results lead to two outlooks. First, a new mold growth vulnerability class could be added to the original mold growth VTT model, namely ``\emph{very very vulnerable}''. Nevertheless, the quality of the parameter estimation is not satisfactory. It is necessary to analyze the model in greater detail.


\section{Reliability of the VTT model to predict the physical phenomenon}
\label{sec:Reliability_VTT_direct}

In the previous section, the robustness of the model was analyzed in order to estimate the vulnerability class of an innovative material using experimental observations. The purpose is now to discuss the reliability of the VTT model.

In order to discuss the reliability, three computations of the VTT model are performed according to \cite{Viitanen2010}, where where a material from the medium resistance class was exposed to constant relative  $\phi \egal 0.97$ and temperature $T \egal 22 \ \mathsf{^{\circ}C} \,$. By slightly changing the value of one parameter ($1\%$ and $0.1\%$) considered as given in the literature, and keeping all the others parameters identical, one looks for the solution $M\,(t)$ with these slightly perturbed parameters.

\begin{figure}
  \centering
  \includegraphics[width=0.45\textwidth]{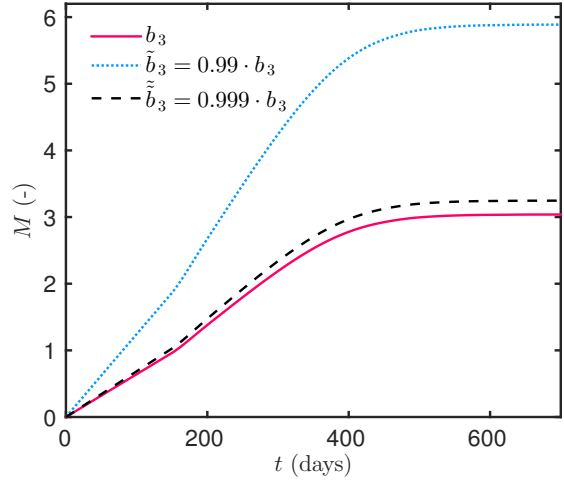}
  \caption{\small\em Comparison of the numerical results from the VTT model.}
  \label{fig:M_fbtilde}
\end{figure}

More precisely, the first one uses the original value $b_{\,3}$ as seen in Eq.~\eqref{eq:function_f}, while the second and the third use a modified parameter $\tilde{b}_{\,3} \egal 0.99 \cdot b_{\,3}$ and $\tilde{\tilde{b}}_{\,3} \egal 0.999 \cdot b_{\,3} \,$, respectively. Figure~\ref{fig:M_fbtilde} shows the variation of the mold growth index $M$ for $700\ \mathsf{days}$. As reported in Table~\ref{tb:res_VTT_fb}, a $1\%$ modification of the parameter leads to a relative error of almost $100\%$ on the results. Furthermore, the two decimals defining the parameters in the VTT model~\eqref{eq:odeM} are not relevant taking into account the precision of the results. These are intrinsic to the VTT model and mainly its function $f\,$. Using the definition of function $f$ in Eq.~\eqref{eq:function_f}, the right-hand side of Eq.~\eqref{eq:odeM} is multiplied by a factor $\exp(\,-\,66.02\,) \ \sim \ 10^{\,-29}\,$, which surely introduces computational rounding errors.

One must acknowledge that the coefficients are never known with less than $1\%$ error. Then the model is so sensitive to the value of its parameters, that its conclusions are not reliable. These observation may lead to the conclusion that the VTT mathematical formulation of mold growth lacks of accuracy.

\begin{table}
\centering
\begin{tabular}{c|c|c|c}
\hline
\hline
\textit{Parameters of function $f$} 
& $b_{\,3}$ 
& $\tilde{b}_{\,3} \egal 0.99 \cdot b_{\,3}$ 
& $\tilde{\tilde{b}}_{\,3} \egal 0.999 \cdot b_{\,3}$  \\
\hline
\hline
\textit{Results $M(\,t \egal 700 \ \mathsf{days}\,)$} 
& $3.03$ 
& $5.88$ 
& $3.24$  \\
\textit{Relative error on $M(\,t \egal 700 \ \mathsf{days}\,)$} 
& -- 
& $94\%$ 
& $7\%$  \\
\textit{$\mathscr{L}_{\,2}$ error} 
& -- 
& $3.3 \cdot 10^{\,4}$ 
& $175$  \\
\hline
\hline
\end{tabular}\bigskip
\caption{\small\em Results of the computation of the VTT model for three values of parameter $b_{\,3}\,$.}
\label{tb:res_VTT_fb}
\end{table}

\section{New mathematical model}
\label{sec:improvement_model}

As highlighted in previous sections, the mathematical formulation of the VTT model in Eq.~\eqref{eq:odeM} may be overly sensitive. As a consequence, the model lacks accuracy when estimating relevant parameters or predicting physical phenomena. Hence, a new mathematical formulation is first proposed to address this issue. Then the robustness of this improved model will be analyzed under both aspects: (i) the estimation of the relevant parameters using the observations from the bamboo fiberboard and (ii) the prediction of mold growth based on data from the literature.


\subsection{Mathematical formulation}

The new formulation proposed is based on the following ordinary differential equation of the mold growth index $M$:
\begin{align}\label{eq:odeM_improved}
  \od{M}{t} \egal k\, (\,T \,,\, \phi \,) \cdot M \cdot \bigl(\, \Minf\, (\,T \,,\, \phi \,) \moins M \,\bigr) \,,
\end{align}
where $\Minf\,$, by analogy with $\Mmax\,$, is the maximum mold growth value for the given temperature and relative humidity conditions, and $k$ is the rate of mold growth. This equation is also known as the logistic equation \cite{Verhulst1845}, particularly used in mathematical biology and population dynamics.

In the case of the VTT model, the values of $M$ computed with this model start from an initial value $M(\,t\egal 0\,) \ > \ 0\,$. However, the improved model increases monotonically to a maximum value $\Minf \,$, where the VTT model assumes a different mold growth rate for $M \ < \ 1$ and for $M \ > \ 1\,$. The parameter $\Minf$ accurately represents the experimental observations of mold growth shown for instance in Figure~\ref{fig:Dataexp_proj2} and typical of the logistic functions. Parameters $\Minf$ and $k$ can include the dependency on the temperature, the relative humidity and other parameters such as the surface quality. Parameter $k$ could be positive or negative or change of signs to represent the physical phenomena of increases or decreases in mold.

Another interesting feature may be noted. For a constant parameter $k\,$, corresponding to fixed temperature and relative humidity conditions, an analytical solution of Eq.~\eqref{eq:odeM_improved} can be readily obtained:
\begin{align}\label{eq:M_solution_analytique}
  M\,(t) \egal \frac{\Minf \, M_{\,0}}{\bigl(\,\Minf \moins M_{\,0} \,\bigr) \cdot \exp \bigl(\, - k \, \Minf \, t\,\bigr) \plus M_{\,0}} \,,
\end{align}
where $M_{\,0} \eqdef M \, (t\egal 0)$ is the initial state of the mold growth, which is a parameter of the improved mathematical model and not from the physical model.

It is important to notice that the model requires a non-zero initial condition $M_{\,0} \ \neq \ 0$ for the process to start. Moreover, the sensitivity of $M$ to the initial condition $M_{\,0}$ is given by:
\begin{align*}
  \od{M}{M_{\,0}} \egal \frac{\Minf^{\,2} \,\exp \bigl(\, - k \, \Minf \, t\,\bigr)  }{%
  \Bigl(\,\bigl(\, \Minf \moins M_{\,0} \,\bigr) \, \exp \bigl(\, - k \, \Minf \, t\,\bigr) \plus M_{\,0} \,\Bigr)^{\,2} } \ .
\end{align*}
It can be noted that the sensitivity tends to $0$ exponentially when $t \ \rightarrow \ \infty \,$.


\subsection{Reliability of the improved model to estimate relevant parameters}

For this study, it should be noted that parameters $k$ and $\Mmax$ are defined as constant for each experiments. In other words, they do not depend on temperature $T$ or relative humidity $\phi\,$. Indeed, it would be possible to have parameters that depend on temperature $T$ and relative humidity $\phi\,$. Consequently, the identification's quality would be improved, since the number of parameters in the model increases for the same amount of measurements. However, the purpose is to discuss the reliability of the improved model. Future works should provide more experiments to work on the variation of these parameters with $T$ and $\phi$ since it is of major importance to have model defined for practical case of building applications. After this remark, the reliability of the improved model must now be evaluated. First, as in the VTT case, it is necessary to project the discrete experimental data before determining the parameter using an optimization process. This is done by projecting the experimental data on continuous functions. However, the choice of the family of functions is based on the analytical solution of Eq.\eqref{eq:M_solution_analytique}.


\subsubsection{Projection of the mold growth data}
\label{sec:proj_mold_data_1}

\begin{figure}
\centering
\subfigure[$\phi_{\,1} \egal 0.75$ \label{fig:exp1_proj2}]{\includegraphics[width=0.45\textwidth]{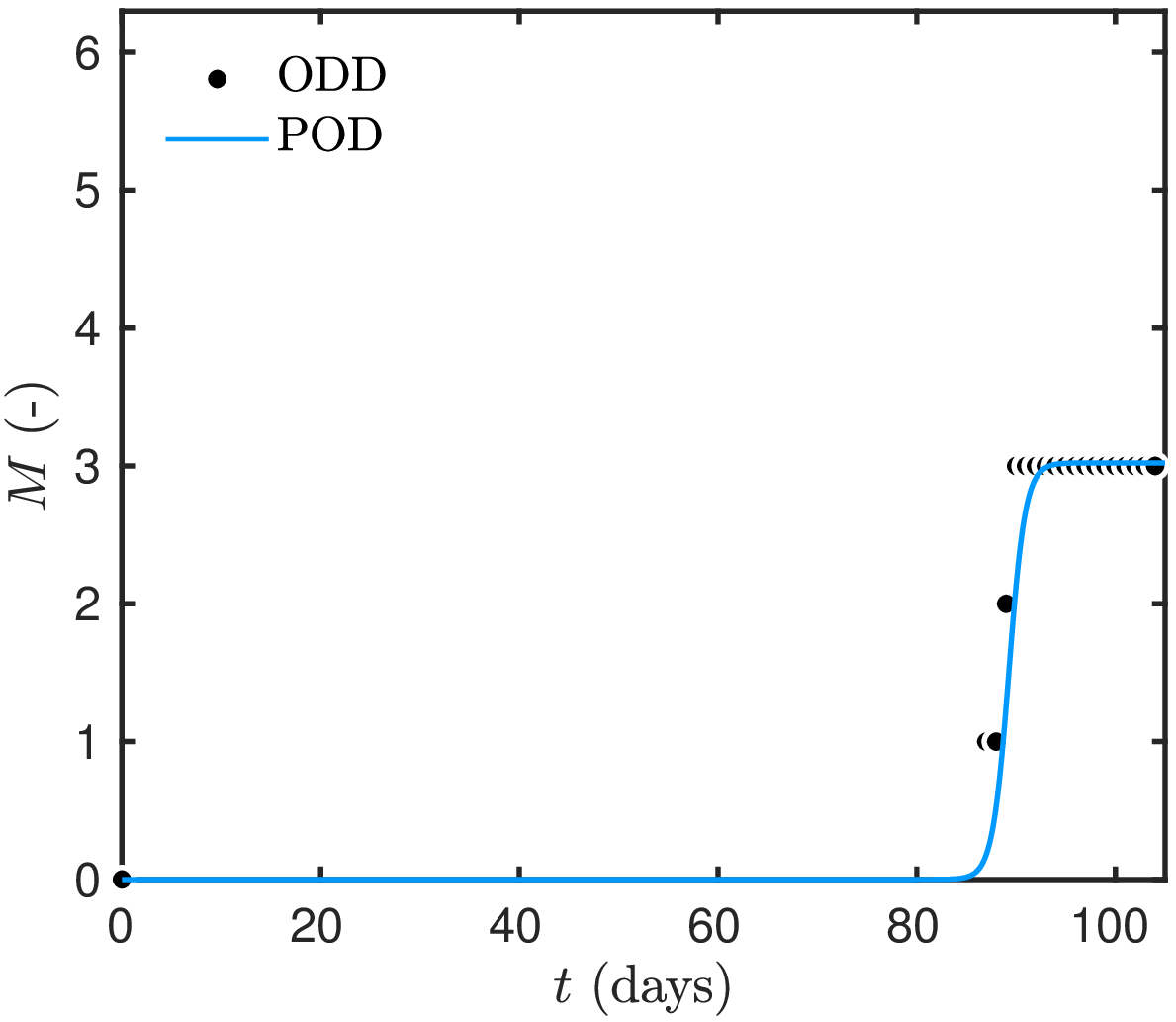}} \hspace{0.2cm}
\subfigure[$\phi_{\,2} \egal 0.84$ \label{fig:exp2_proj2}]{\includegraphics[width=0.45\textwidth]{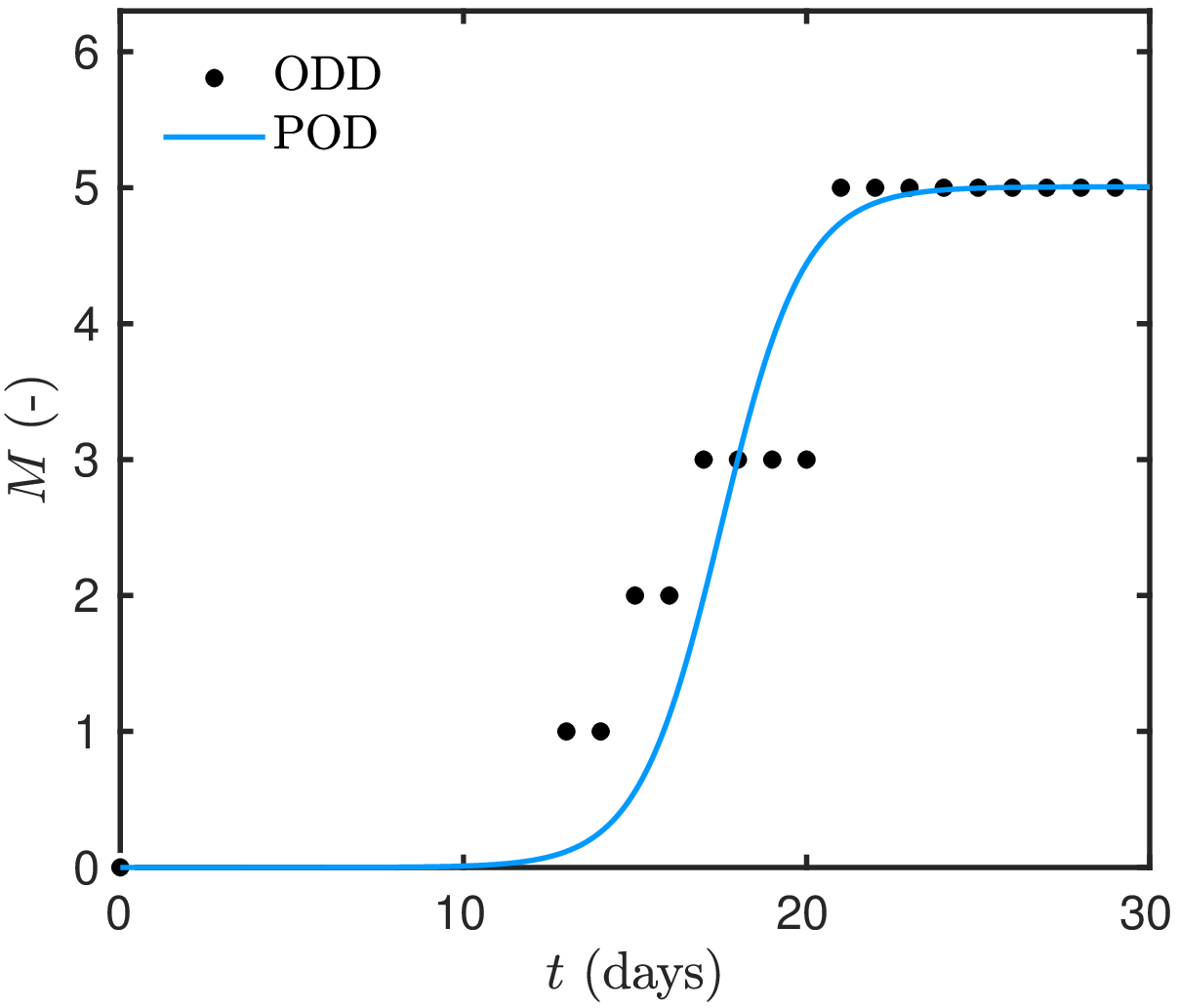}} \\
\subfigure[$\phi_{\,3} \egal 0.97$ \label{fig:exp3_proj2}]{\includegraphics[width=0.45\textwidth]{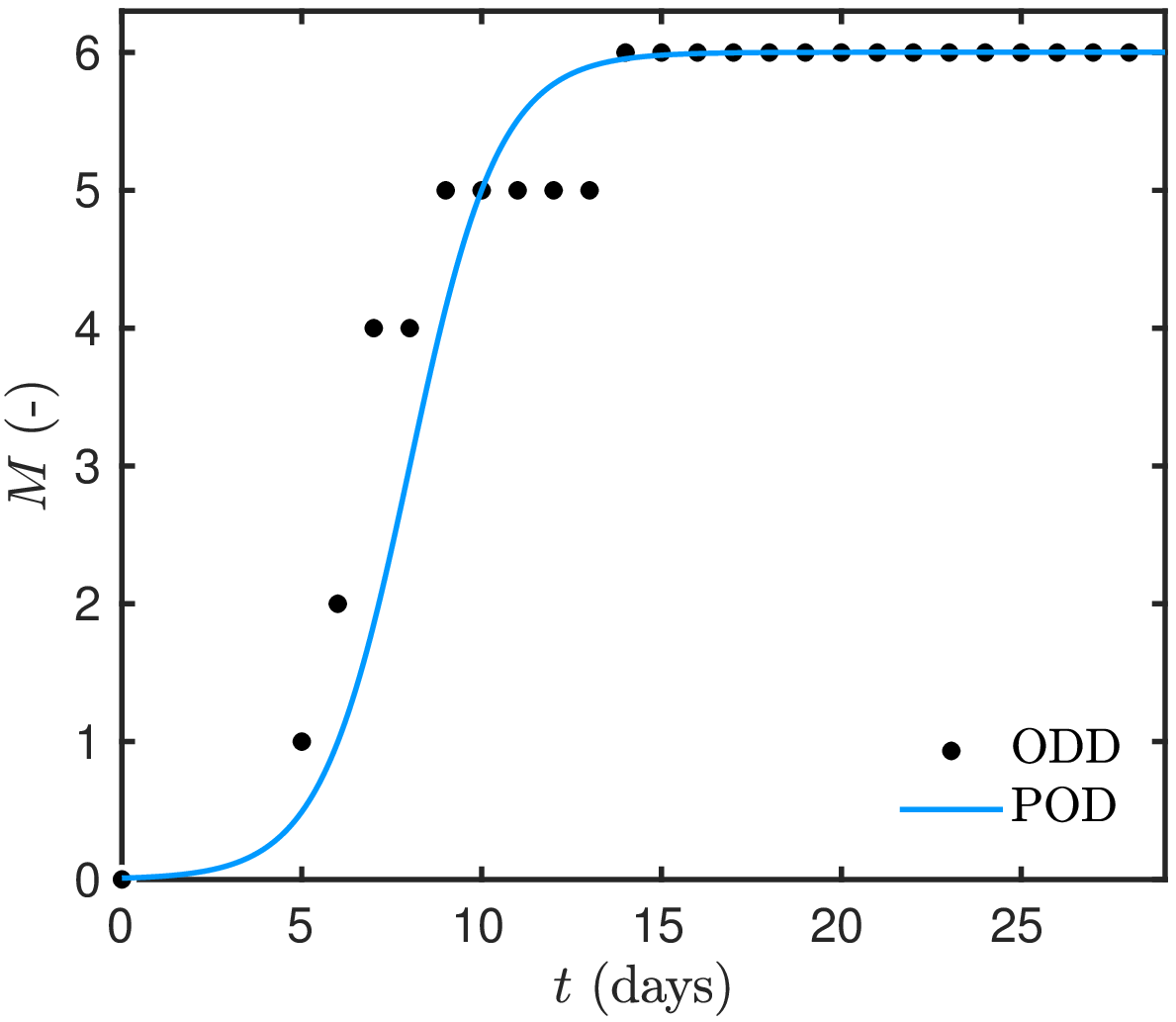}}
\caption{\small\em Experimental ODD available for three different values of relative humidity and the POD resulting from the projection.}
\label{fig:Dataexp_proj2}
\end{figure}

As for the VTT model, the continuous model from Eq.~\eqref{eq:odeM_improved} cannot match the discrete experimental data. Thus, the experimental data are projected using the logistic function:
\begin{align}\label{eq:projection2}
  h_{\,a,\,b,\,c\,}(\,t\,) \egal \frac{a}{1 \plus \exp \bigl(\, - \, b \cdot (\,t \moins c \,) \,\bigr)} \,.
\end{align}
The values of parameters $\bigl(\,a\,,\,b\,,\,c\,\bigr)$ are given in Table~\ref{tb:coeff_proj2}. It can be noted that the family functions used for the projection are different from the one chosen in Section~\ref{sec:projected_mold_growth_data}. The choice of the projection functions is crucial since it is important not to introduce errors while performing the projection of the ODD. Moreover, it has to be chosen in accordance with the type of solutions used in the mathematical model. The time domain is discretized with a uniform mesh using a time step of $6 \ \mathsf{min}\,$.  Figures~\ref{fig:Dataexp_proj2}(a,b,c) shows the comparison between the projected and the discrete experimental data.

\begin{table}
\centering
\begin{tabular}{c|ccc|cc}
\hline
\hline
\textit{Parameters} 
& \textit{Exp.} $\phi_{\,1}$ 
& \textit{Exp.} $\phi_{\,2}$ 
& \textit{Exp.} $\phi_{\,3}$
& Data from \cite{Johansson2013}
& Data from \cite{Nielsen2004}\\
\hline
\hline
$a$ & $3.02$ & $5.01$ & $6.00$ & $5.24$ & $2.58$\\
$b$ & $1.21$ & $0.83$ & $0.81$ & $0.48$ & $0.53$\\
$c$ & $89.27$ & $17.50$ & $7.98$ & $17.26$ & $14.84$\\
\hline
\hline
\textit{$\L_{\,2}$ residual}
& $0.52$ & $0.60$ & $0.56$ & $0.2$ & $0.4$ \\
\hline
\hline
\end{tabular}\bigskip
\caption{\small\em Parameters for the projection of the experimental data using the family of functions $h_{\,a,\,b,\,c\,}(\,t\,)$ to obtain the projected observed data for Section~\ref{sec:proj_mold_data_1} ($\phi_{\,1}\,,\,\phi_{\,2}\,,\,\phi_{\,3}$) and Section~\ref{sec:proj_mold_data_2} (data from \cite{Johansson2013, Nielsen2004}).}
\label{tb:coeff_proj2}
\end{table}


\subsubsection{Parameter estimation}

Using the projected data, the purpose is to estimate the parameters $\bigl(\, \Minf \,,\, k \, \bigr)\,$. According to its definition, the parameter $\Minf$ can be directly determined since it is the maximum mold growth of the observed discrete data. On the other hand, it is interesting to note that from Eq.~\eqref{eq:M_solution_analytique}, an analytical expression of the time $t_{\,1}$ to reach $M \egal 1$ can be computed:
\begin{align*}
  t_{\,1} \ \eqdef \arg_{\,t} \bigl(\, \,M\,(\,t\,) \egal 1 \,\bigr) \egal \frac{1}{k \, \Minf} \; \ln \Biggl(\,\frac{M_{\,0} \bigl(\, \Minf \moins 1 \,\bigr)}{\Minf \moins M_{\,0}}\,\Biggr)\,.
\end{align*}
Using this expression, it is possible to estimate parameter $k$ directly. However, as mentioned above, the model is sensitive to the initial condition. Most particularly, the derivative of $t_{\,1}$ relative to $M_{\,0}$ is:
\begin{align*}
  \od{t_{\,1}}{M_{\,0}} \egal \frac{1}{k \, M_{\,0} \, \bigl(\, M_{\,0} \moins \Minf \,\bigr)} \,.
\end{align*} 
Since $M_{\,0} \ \ll \ 1$, the derivative is large $\od{t_{\,1}}{M_{\,0}} \ \gg \ 1 \,$. In addition, the experimental data may lack accuracy when estimating $t_{\,1}\,$. To circumvent this problem, the least square method is used to estimate the parameters $\bigl(\, \Minf \,,\, k \, \bigr)\,$. To do this, a cost function is defined according to Eq.~\eqref{eq:cost_function} and minimized using the \texttt{fmincon} algorithm from the \texttt{Optimization Toolbox} in the \texttt{Matlab\texttrademark} environment.


\subsubsection{Results and discussion}

\begin{figure}
\centering
\subfigure[$\phi_{\,1} \egal 0.75$ \label{fig:M1_ft_IP2}]{\includegraphics[width=0.45\textwidth]{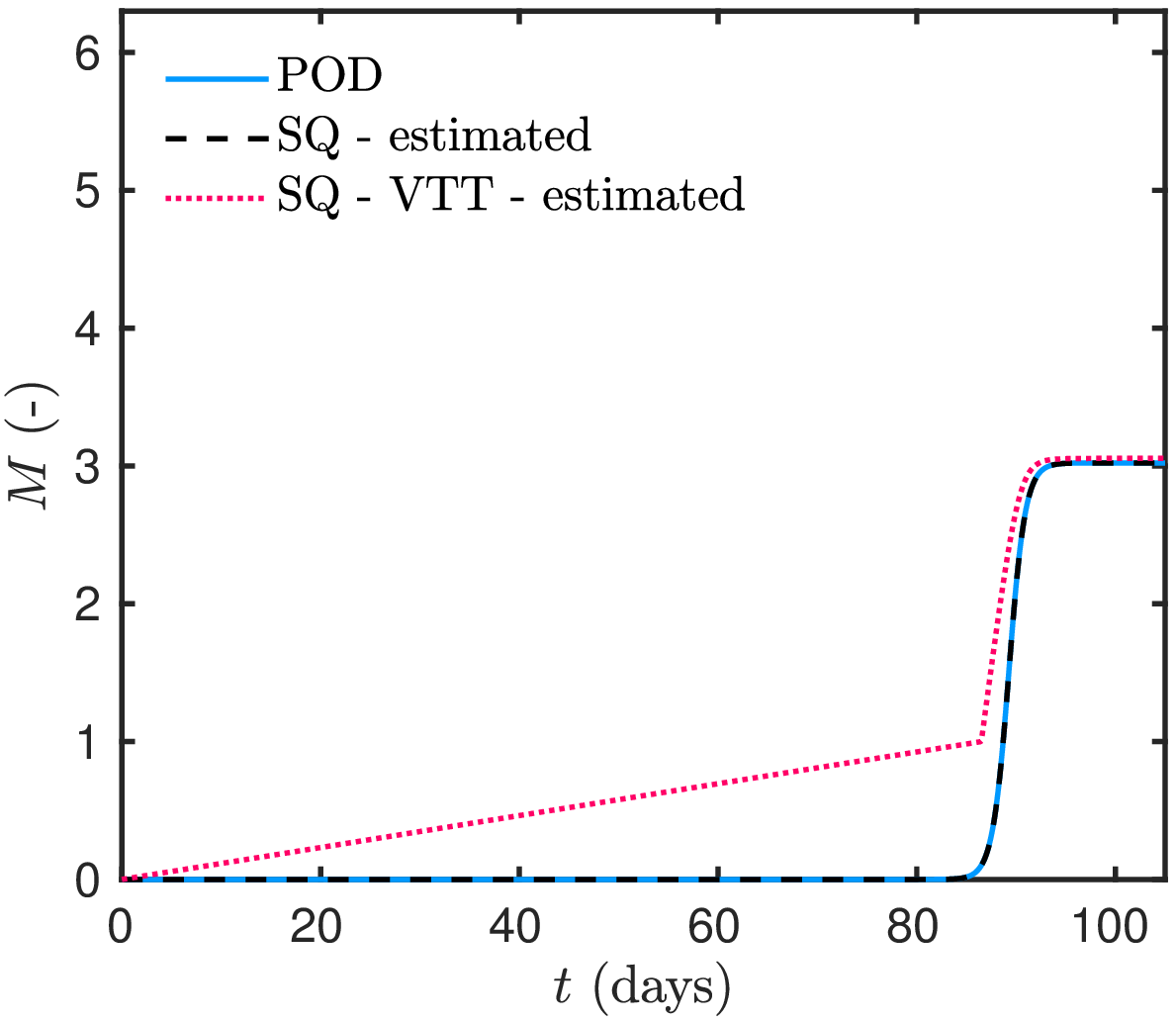}} \hspace{0.2cm}
\subfigure[$\phi_{\,2} \egal 0.84$ \label{fig:M2_ft_IP2}]{\includegraphics[width=0.45\textwidth]{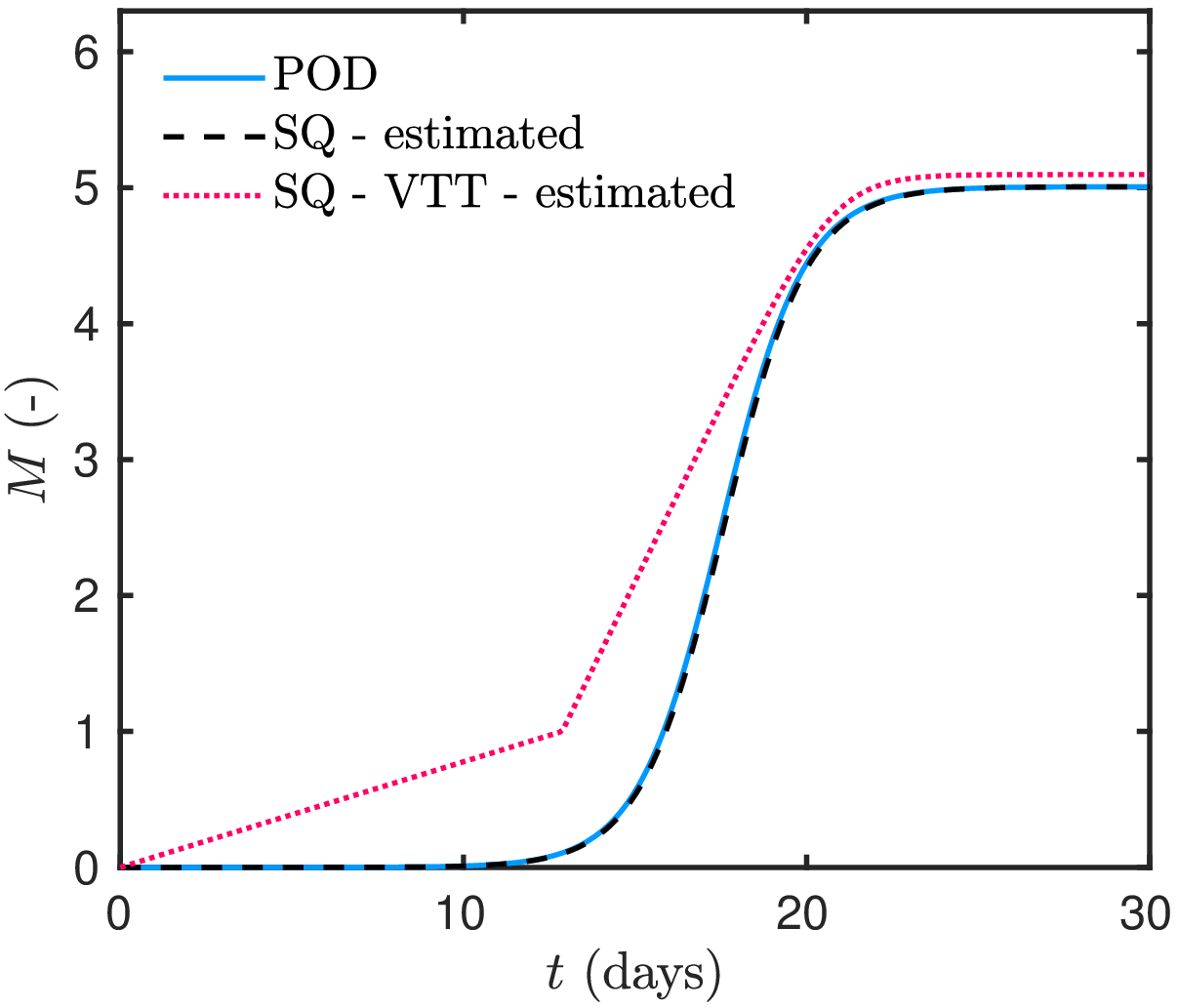}} \\
\subfigure[$\phi_{\,3} \egal 0.97$ \label{fig:M3_ft_IP2}]{\includegraphics[width=0.45\textwidth]{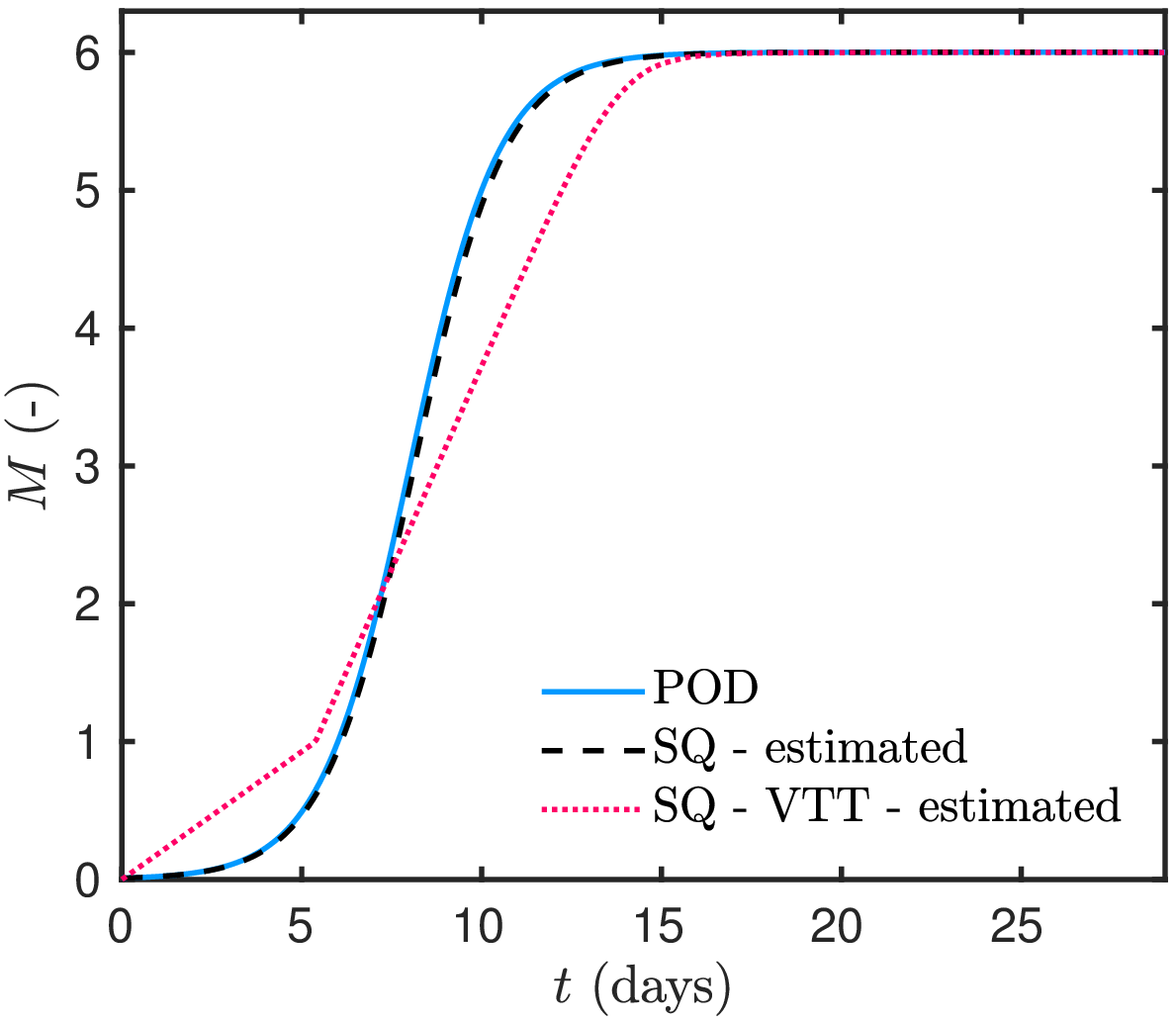}}\hspace{0.2cm}
\subfigure[\label{fig:res_ft_IP2}]{\includegraphics[width=0.45\textwidth]{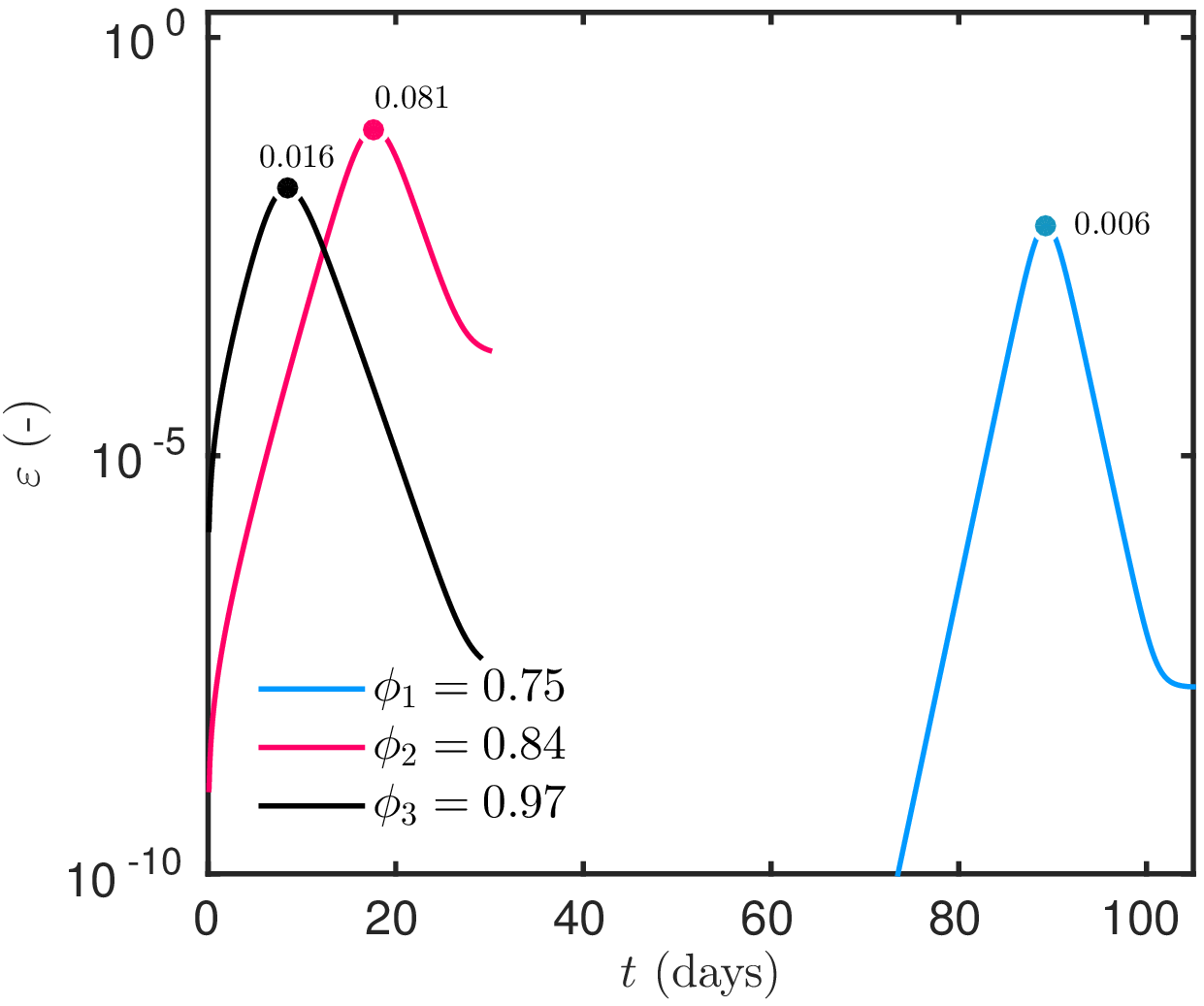}} 
\caption{\small\em \emph{(a,b,c)} Comparison of the POD with the SQ using the improved mold growth model and the VTT model with the estimated parameters. \emph{(d)} Residual between the POD data and the simulated ones using the estimated parameters.}
\label{fig:M_ft_IP2}
\end{figure}

In terms of parameter estimation, where the model in Eq.~\eqref{eq:odeM} required the determination of five parameters, only three are needed for the improved model: $\Minf\,$, $k\,$ and $\Mzero\,$. The parameter $\Minf$ is directly obtained from the maximum value of $M$ observed for a fixed temperature and a given relative humidity when the steady state is reached. The other two parameters $k$ and $\Mzero$ are estimated by minimizing the norm of the differences between the projected observed data and the simulated quantities, defined in Eq.~\eqref{eq:cost_function}. The only constraint solely concerns $\Minf$ according to its definition: $\Minf \ \leqslant \ 6\,$.
The estimated parameters $\Mzero^{\,\circ}\,$, $k^{\,\circ}$ and $\Minf^{\,\circ}$ are reported in Table~\ref{tb:unknown_parameters_IP2}. Compared to the VTT model parameters, $k^{\,\circ}$ and $\Minf^{\,\circ}$ have the same order of magnitude and increase with the relative humidity. Figure~\ref{fig:M_ft_IP2}(a,b,c) compares the POD with the SQ obtained with the estimated parameters. Very good agreement can be noted. The residual, shown in Figure~\ref{fig:res_ft_IP2}, is lower than for the previous estimation. The maximum value of the residual is of the order $\O(\,10^{\,-2}\,)$. The $\L_{\,2}$ norm of the residual is $8 \cdot 10^{\,-4}\,$ $1.4 \cdot 10^{\,-3}$ and $2.1 \cdot 10^{\,-3}$ for each experiment, which is much lower than for identification for the VTT model. Moreover, the results of the original VTT model with its estimated parameters are recalled in Figure~\ref{fig:M_ft_IP2}(a,b,c).  The improved model provides better agreement with the projected observed data than the original VTT.

As reported in Table~\ref{tb:unknown_parameters_IP2}, the estimation error is much lower for this model than for the original VTT model. Here the initial condition was considered as an unknown parameter. One of the drawbacks of the improved model is its sensitivity to the initial value of the problem. This could be solved fixing the initial value for all types of materials or by using no exponential in the projection of the ODD nor in the evolution of the mold growth model.

\begin{table}
\centering
\begin{tabular}{c|c|c|c}
\hline
\hline
\textit{Experiment} & $\Mzero^{\,\circ}$ & $k^{\,\circ}$ & $\Minf^{\,\circ}$ \\[3pt]
\hline
\hline
$\phi_{\,1} \egal 0.75$ & $4 \cdot 10^{\,-7}$ & $0.39$ & $3.02$ \\
\textit{Error estimator} $\eta$
& $1.4 \cdot 10^{\,-12}$
& $3.2 \cdot 10^{-3}$
& $1.2 \cdot 10^{-2}$ \\
\hline
$\phi_{\,2} \egal 0.84$ & $2 \cdot 10^{\,-6}$ & $0.16$ & $5.00$ \\
\textit{Error estimator} $\eta$
& $1.6 \cdot 10^{\,-7}$
& $7.5 \cdot 10^{-3}$
& $3.5 \cdot 10^{-3}$ \\
\hline
$\phi_{\,3} \egal 0.97$ & $9 \cdot 10^{\,-3}$ & $0.13$ & $6.00$ \\
\textit{Error estimator} $\eta$
& $2.1 \cdot 10^{\,-3}$
& $4.6 \cdot 10^{-3}$
& $2.5 \cdot 10^{-3}$ \\
\hline
\hline
\end{tabular}\bigskip
\caption{\small\em Estimated parameters of the improved model for the bamboo fiberboard.}
\label{tb:unknown_parameters_IP2}
\end{table}


\subsection{Reliability of the improved model to predict the physical phenomenon}

First, the reliability of the improved model is evaluated using data from the literature. Since this analysis concerns only one type of material (pine wood) and many other materials can be used for building applications, the possibility of defining vulnerability classes of material is studied in the next Section.


\subsubsection{Using data from the literature}
\label{sec:proj_mold_data_2}

\begin{figure}
  \centering
  observed data \cite{Johansson2013} 
  \hspace{4.cm}
  observed data \cite{Nielsen2004} \\
  \subfigure[\label{fig:MJOH_ft}]{\includegraphics[width=0.43\textwidth]{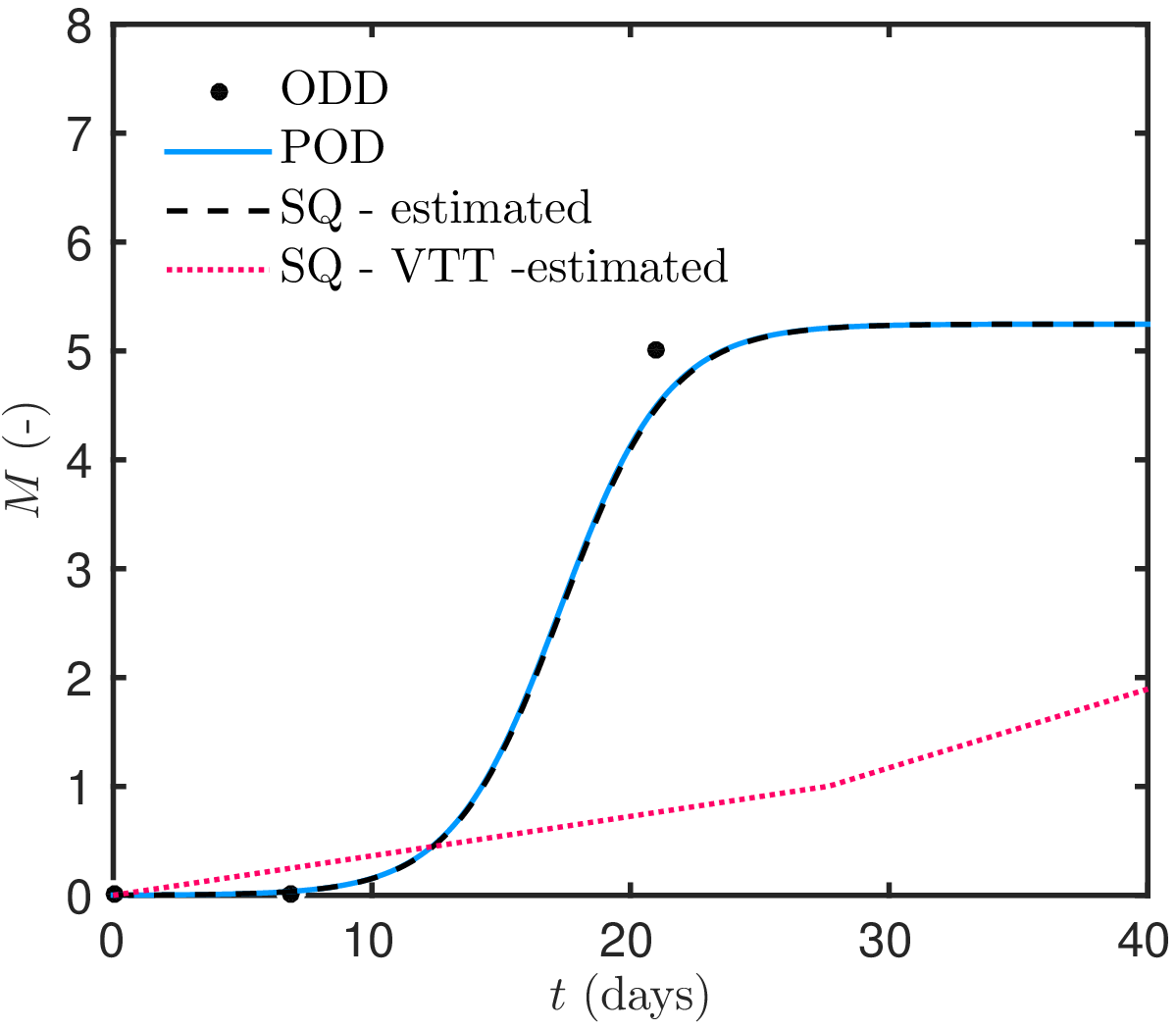}} \hspace{0.2cm}
  \subfigure[\label{fig:MNIELSEN_ft}]{\includegraphics[width=0.45\textwidth]{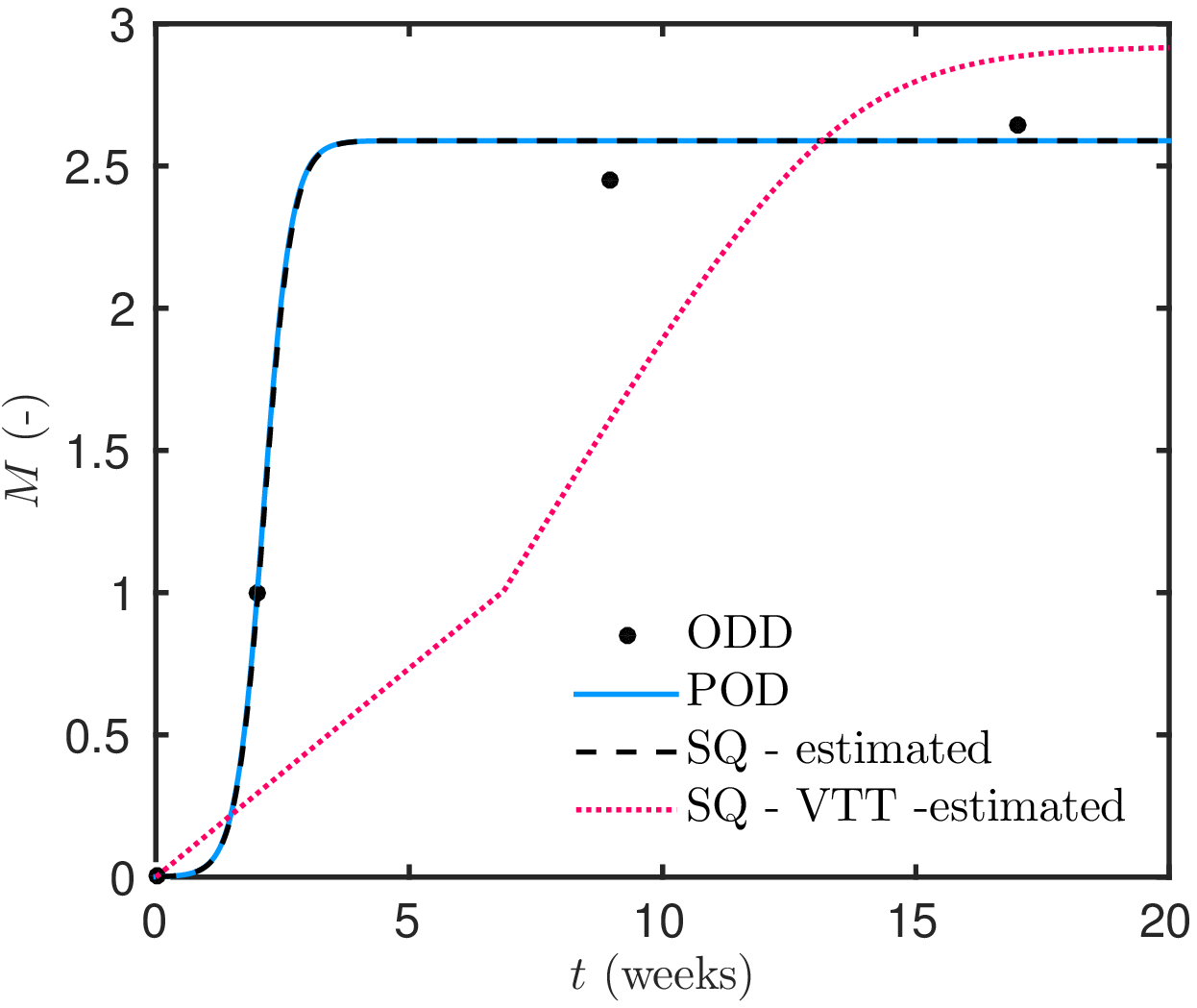}}\\
  \subfigure[\label{fig:resJOH_ft}]{\includegraphics[width=0.45\textwidth]{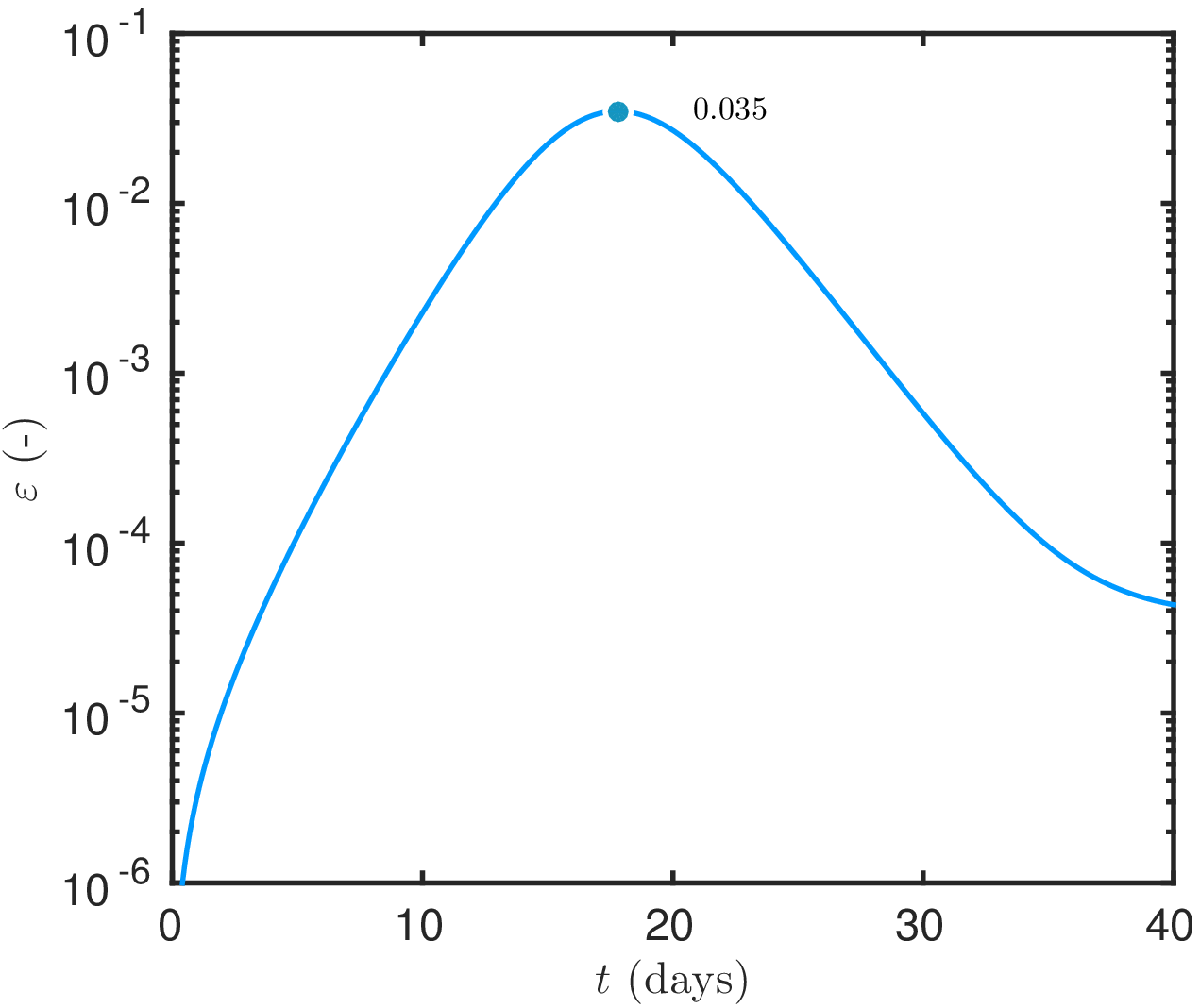}} \hspace{0.2cm}
  \subfigure[\label{fig:resNIELSEN_ft}]{\includegraphics[width=0.45\textwidth]{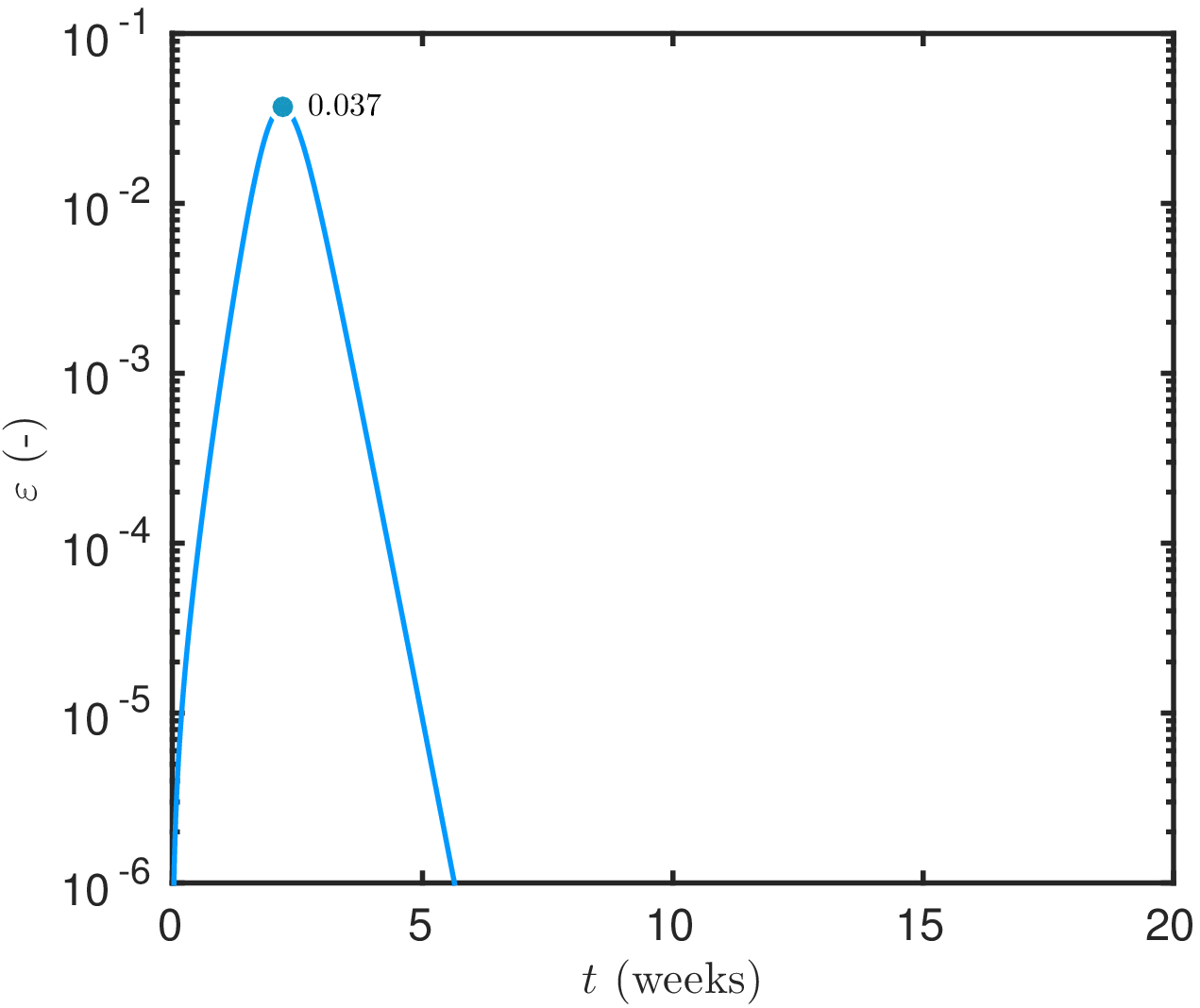}} 
  \caption{\small\em \emph{(a,b)} Comparison of the POD with the ODD and the SQ using the improved mold growth model with the estimated parameters. \emph{(c,d)} Residual between the POD and the SQ with the estimated parameters. The observed data are taken from \cite{Johansson2013} \emph{(a,c)} and \cite{Nielsen2004} \emph{(b,d)}.}
  \label{fig:M_ft_IP3}
\end{figure}

To evaluate the reliability of the improved model to predict the physical phenomena, two sets of experimental data from the literature are used to estimate the unknown parameters $\bigl(\,k \,,\, \Mzero \,,\, \Minf \,\bigr)$ for the improved model. The first set is taken from \cite{Johansson2013}, where a planted pine sapwood is exposed to constant conditions at $T \egal 22\ \mathsf{^\circ C}$ and $\phi \egal 0.9 \,$. The second set comes from \cite{Nielsen2004} for a similar material and constant conditions at $T \egal 25\ \mathsf{^\circ C}$ and $\phi \egal 0.86 \,$. The ODD are shown in Figure~\ref{fig:M_ft_IP3}(a,b). As detailed in previous sections, the discrete data are projected using the function defined in Eq.~\eqref{eq:projection2} and the coefficients are reported in Table~\ref{tb:coeff_proj2}.

The solutions of the parameter estimation problem are reported in Table~\ref{tb:unknown_parameters_IP3}. Again, the estimated parameters have the same order of magnitude. Figure~\ref{fig:M_ft_IP3}(a,b) compares the SQ with the estimated parameter and the POD for \cite{Johansson2013, Nielsen2004}, respectively. For both experiments, satisfactory agreement is noted with a residual of the order $\O\,(\,10^{\,-2}\,)$, as shown in Figures~\ref{fig:M_ft_IP3}(c,d). It can be noted that the original VTT model does not succeed in representing the physical phenomena. These observations agree with those reported in \cite{Vereecken2015} for the same experiments.

\begin{figure}
\centering
\subfigure[\label{fig:MJOH_ft_fM0}]{\includegraphics[width=0.43\textwidth]{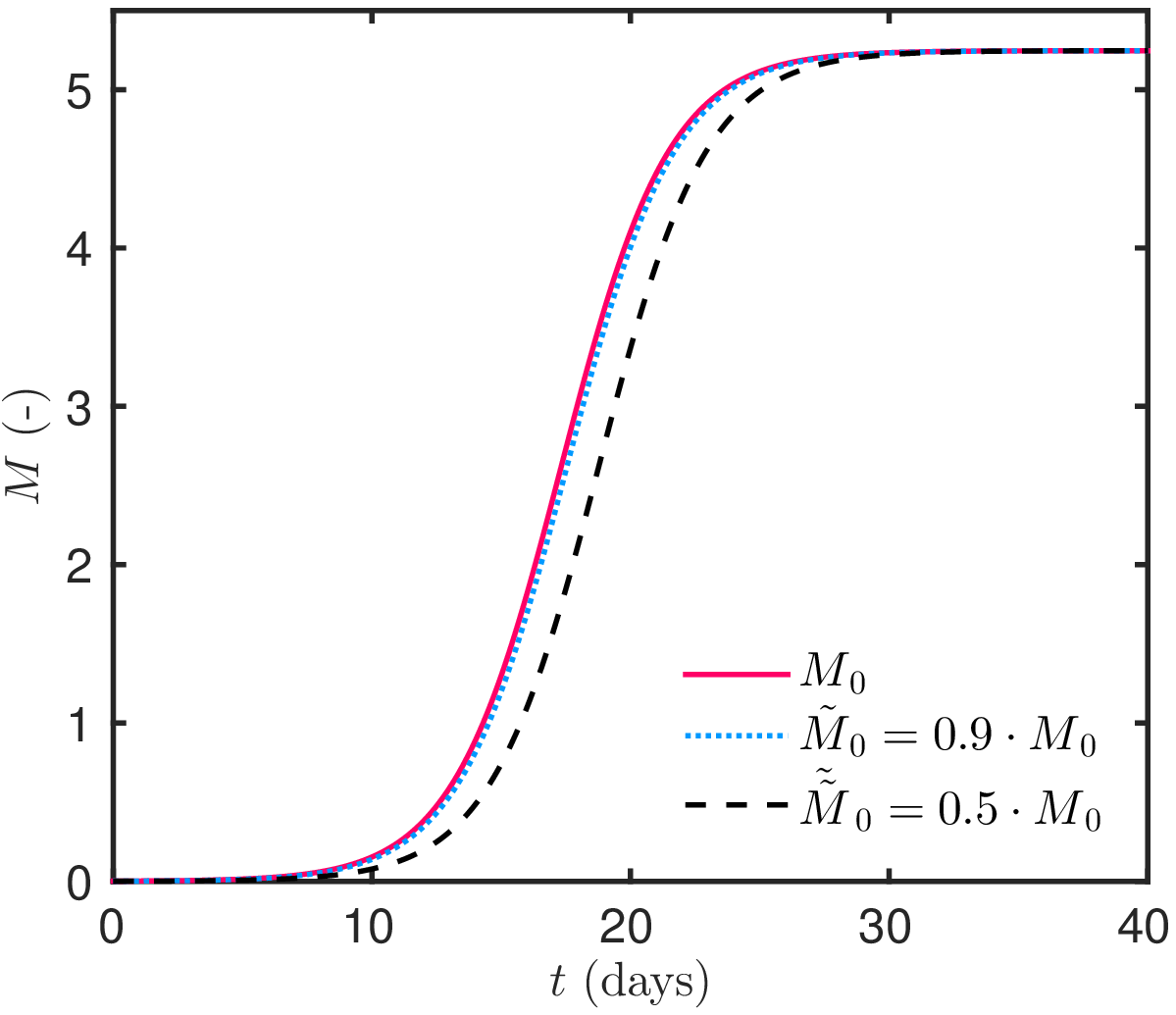}} \hspace{0.2cm}
\subfigure[\label{fig:JOH_err_fM0}]{\includegraphics[width=0.45\textwidth]{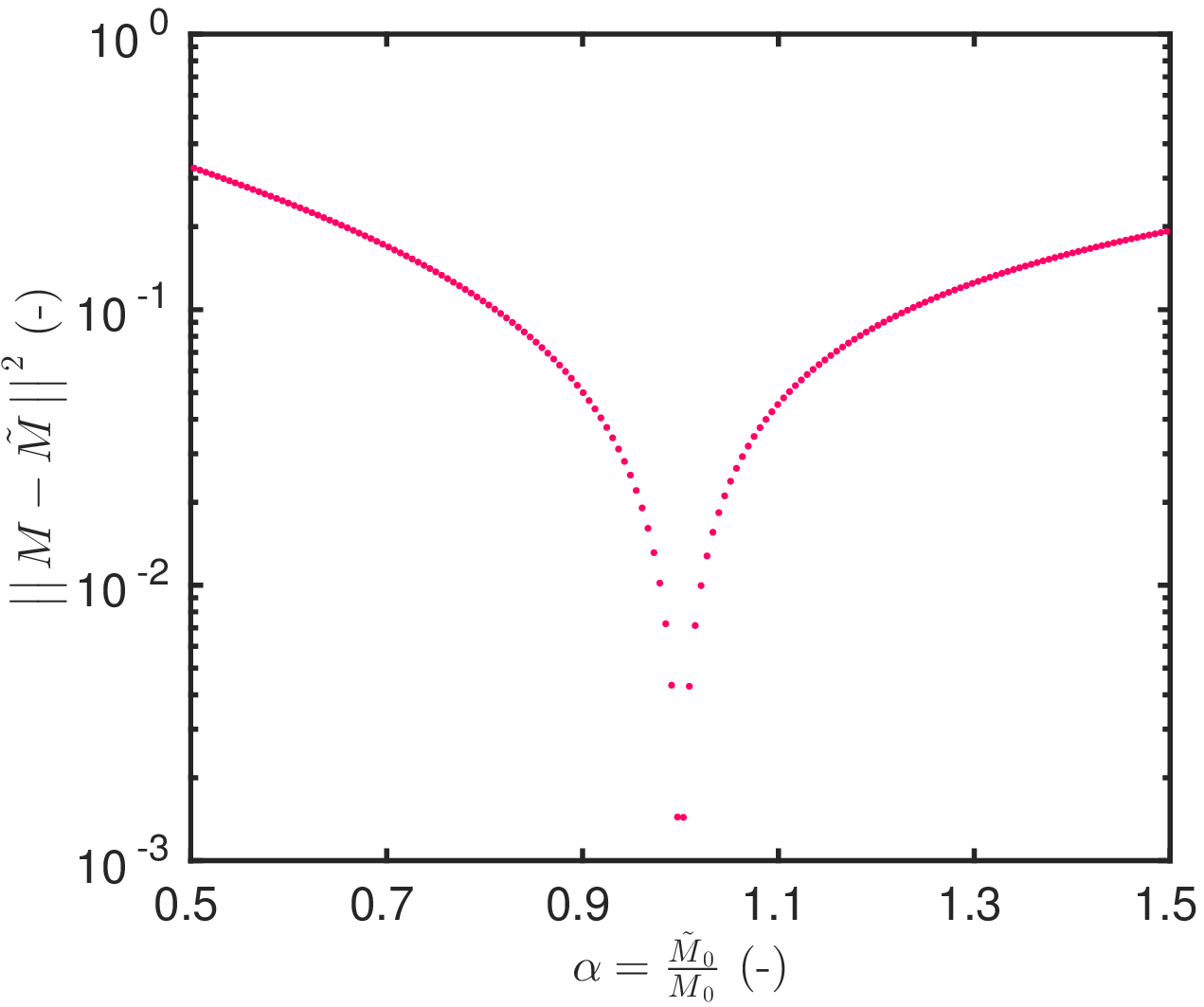}}
\caption{\small\em Comparison of the numerical results from the improved model for different values of the initial condition $\tilde{M}_{\,0} \egal \alpha \cdot M_{\,0}\,$, other parameters remaining constant, for the experiment defined by \cite{Johansson2013}.}
\label{fig:M_fM0}
\end{figure}

Since the improved model is sensitive to the initial condition $M_{\,0}\,$, a parametric study is carried out by performing computations of $M$ for different values of the initial condition $\tilde{M}_{\,0} \egal \alpha \cdot M_{\,0}\,$. These computations are realized for the experiments from \cite{Johansson2013} and all other estimated parameters ($k$ and $\Minf$) remain unchanged. The results are shown in Figure~\ref{fig:M_fM0}(a,b). A $50\%$ modification of the numerical value of the initial conditions does not imply substantial change in the prediction of the physical phenomena. As expected, the initial condition has no influence on the maximum mold growth $\Minf\,$. Figure~\ref{fig:JOH_err_fM0} shows the variation of the $\mathscr{L}_{\,2}$ error with the modification. The magnitude of the error on the prediction is small (at most $20\%$ for a $50\%$ change) compared to that observed ($94\%\,$, $7\%$) by modifying the numerical value of parameter $b_{\,3}$ in the VTT model by $1\%\,$, $0.1\%$ (see Table~\ref{tb:res_VTT_fb}). This explains why that this improved model is more reliable.

These results highlight that the improved mold growth model~\eqref{eq:odeM_improved} can be used to fit the unknown parameters using projected data obtained from other discrete experimental data. It can be noted that, for the bamboo fiberboard and a fixed temperature, parameters $k$ and $\Minf$ decrease and increase with relative humidity, respectively. Further experimental data are required to study the variation of these parameters with temperature and relative humidity.

\begin{table}
\centering
\begin{tabular}{c|c|c|c|c}
\hline
\hline
\textit{Experiment} & \textit{Reference} & $\Mzero^{\,\circ}$ & $k^{\,\circ}$ & $\Minf^{\,\circ}$ \\[3pt]
\hline
\hline
$T \egal 22\ \mathsf{^\circ C}\,$, $\phi \egal 0.9$ & \cite{Johansson2013} & $1.3 \cdot 10^{\,-3}$ & $0.10$ & $5.24$ \\
$T \egal 25\ \mathsf{^\circ C}\,$, $\phi \egal 0.86$ & \cite{Nielsen2004} & $1.0 \cdot 10^{\,-3}$ & $0.19$ & $2.58$ \\
\hline
\hline
\textit{Error estimator} $\eta$
& - 
& $1.2 \cdot 10^{\,-3}$
& $3.4 \cdot 10^{-3}$
& $2.7 \cdot 10^{-3}$ \\[3pt]
\hline
\hline
\end{tabular}\bigskip
\caption{\small\em Estimated parameters of the improved model for pine sapwood using the experimental data from \cite{Johansson2013} and \cite{Nielsen2004}.}
\label{tb:unknown_parameters_IP3}
\end{table}

\subsubsection{On the definition of vulnerability classes}
\label{sec:def_vuln_classes}

\begin{figure}
  \centering
  \subfigure[\label{fig:M_VR_R}]{\includegraphics[width=0.45\textwidth]{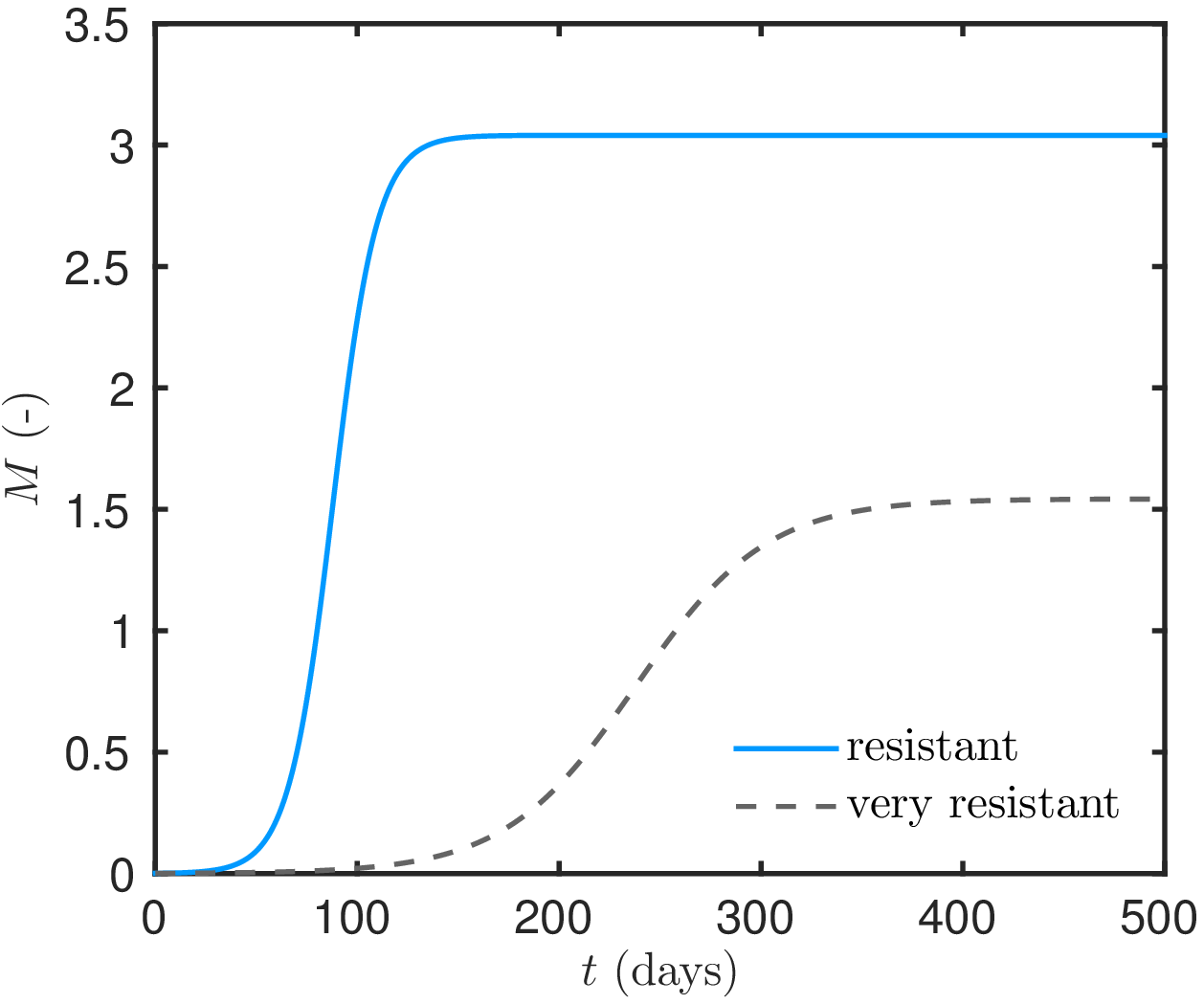}} \hspace{0.2cm}
  \subfigure[\label{fig:M_VS_S}]{\includegraphics[width=0.43\textwidth]{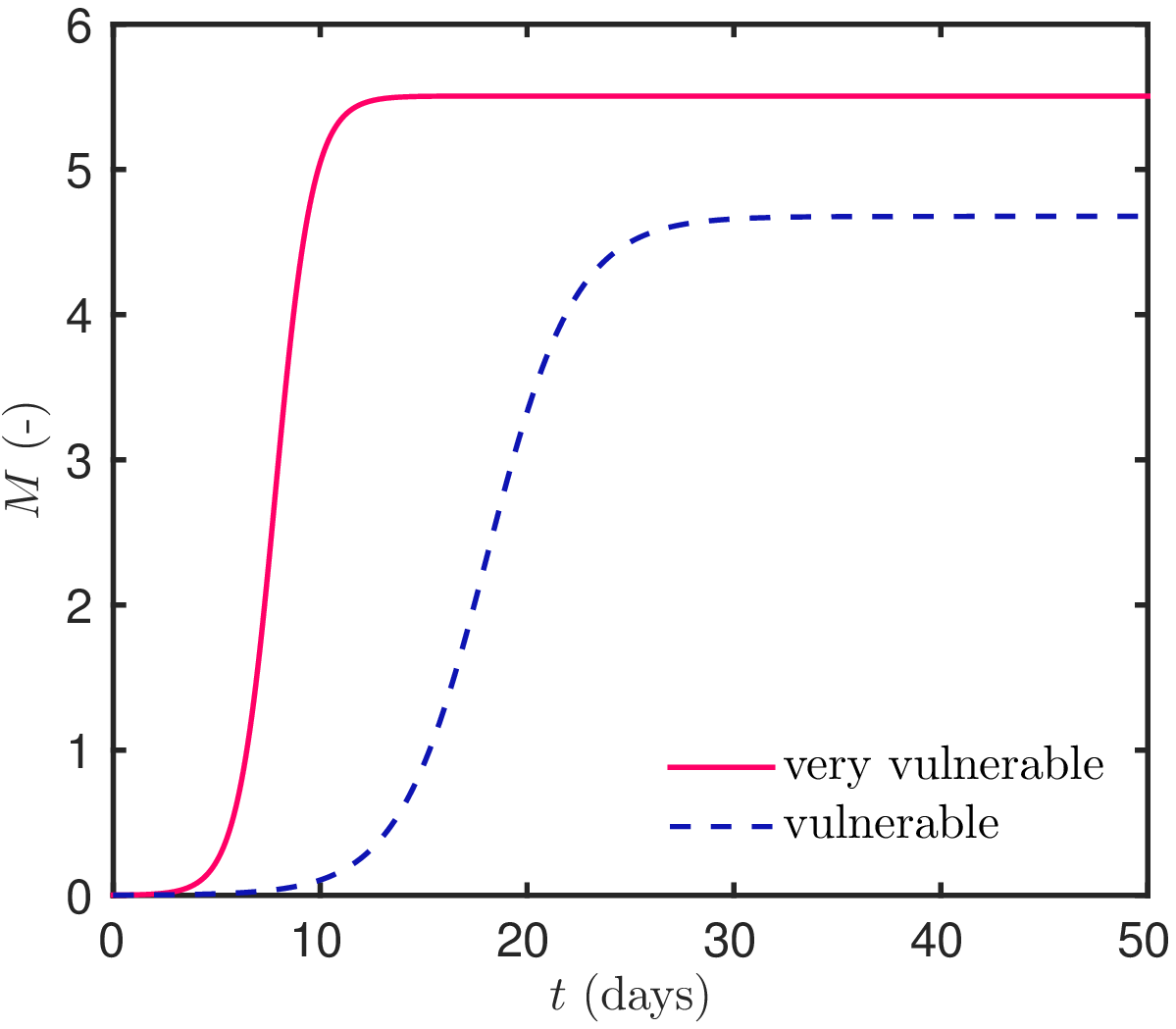}}
  \caption{\small\em Numerical predictions of the improved model for different values of parameters $k$ and $\Minf\,$, reported in Table~\ref{tb:improvedVTT_fclasses}.}
  \label{fig:improvedVTT_fclasses}
\end{figure}

The previous section highlighted that the improved model accurately predicts the physical phenomenon based on experimental data from the literature. This study focuses on a single type of material and one can wonder if vulnerability classes could be defined for the improved model. To answer this question, several computations of the improved model are done considering a fixed initial value $\Mzero \egal 10^{\,-3}\,$. As reported in Table~\ref{tb:improvedVTT_fclasses}, several values of parameters $k$ and $\Minf$ are considered for the various vulnerability classes. These parameters decrease when the material is less vulnerable to mold growth. It should be noted that these parameters are valid for fixed conditions, particularly for temperature and relative humidity. Figure~\ref{fig:improvedVTT_fclasses}(a,b) shows the time variation of $M$ for the different vulnerability classes. Mold growth increases faster for the more vulnerable materials and reaches a higher level of $\Minf \,$. For a \emph{very vulnerable} material, mold growth starts after $5 \ \mathsf{days}$ of exposure while for a \emph{very resistant} one, the process appears after $200 \ \mathsf{days}\,$. These results highlight the possibility of defining mold growth vulnerability classes for materials with the improved model.

\begin{table}
  \centering
  \begin{tabular}{l|c|c}
  \hline
  \hline
  \textit{Vulnerability classes} & $k^{\,\circ}$ & $\Minf^{\,\circ}$ \\[3pt]
  \hline 
  \hline
  \textrm{Very vulnerable} & $0.2$ & $5.51$ \\
  \textrm{Vulnerable} & $0.1$ & $4.67$ \\
  \textrm{Resistant} & $0.03$ & $3.04$ \\
  \textrm{Very resistant} & $0.02$ & $1.54$ \\
  \hline
  \hline
  \end{tabular}\bigskip
  \caption{\small\em Parameters of the improved model as a function of the vulnerability class.}
  \label{tb:improvedVTT_fclasses}
\end{table}


\section{Conclusion}
\label{sec:final_remarks}

In buildings, excessive levels of moisture may lead to several disorders that damage the quality of a construction and deteriorate the occupants' comfort. Mold risks are of capital importance since mold can be toxic for the occupants. Several mold growth models have been proposed in the literature that can predict risks depending on relative humidity and temperature conditions. Among them, the VTT model was proposed for wood and other building materials \cite{Hukka1999, Viitanen2010}.

The reliability of a model can be tested on both its ability to represent/interpolate the experimental data and both its the sensitivity/robustness of its identified parameters computed from the experimental data, if they are identifiable. The reliability of the VTT model is investigated for coditions that are constant over time for temperature and relative humidity in this article.

First of all, in Section~\ref{sec:experimental_data}, three experiments are carried out for a recently developed material, measuring mold growth for a fixed temperature and three levels of relative humidity. No DOD may fit solutions of the Ordinary Differential Equation. So the ODD are projected on intermediate continuous functions (Section~\ref{sec:projection_ODD_VTT}).

Then in Section~\ref{sec:estimation_parameters_VTT} the parameters of the VTT model are identified from the intermediate projected experimental data to estimate five parameters defining a material property denoted as the mold vulnerability class in the original model. The results highlight that the parameters determined are totally out of the range of the vulnerability class defined in the original VTT model. Moreover, the quality of the parameter estimation is not satisfactory.

In Section~\ref{sec:Reliability_VTT_direct}, it is highlighted that a simulation using the VTT model is much too sensitive to one of the parameters of $f$ that is considered as reliable in the literature. These investigations lead to the conclusion that the mathematical formulation of the physical model of mold growth lacks accuracy.

Therefore, in Section~\ref{sec:improvement_model} an improved mold growth model with a new mathematical formulation is studied. It is based on an ordinary differential equation, also called the logistic equation. Three parameters are involved in the formulation: (i) $k\,(\,T\,,\,\phi\,)$ the rate of mold growth, (ii) $\Minf\,(\,T\,,\,\phi\,)$ the maximum mold growth value for the given temperature and relative humidity conditions, and (iii) $\Mzero$ the initial mold growth value. These parameters are estimated from our experimental data on bamboo fiberboard. The parameter estimation quality is very satisfactory. The estimated parameters have the same order of magnitude around unity for the three experiments. The reliability of this model is also tested for two other sets of experimental data from studies reported in the literature. It is shown that the improved formulation of the mold growth model accurately predicts the physical phenomenon based on experiments from other authors. In addition, the possibility of defining mold growth vulnerability classes of materials for the improved model is demonstrated in Section~\ref{sec:def_vuln_classes}.

Further research may focus on producing and using further experimental data to define the functions describing the variation of the parameters $k$ and $\Minf$ with large quantities such as the temperature or relative humidity. An important issue concerns the accuracy of the model to predict mold increase and decrease during transient condition of hygrothermal fields.


\subsection*{Acknowledgments}
\addcontentsline{toc}{subsection}{Acknowledgments}

The authors acknowledge the Junior Chair Research program ``\textit{Building performance assessment, evaluation and enhancement}'' from the University \textsc{Savoie Mont Blanc} in collaboration with The French Atomic and Alternative Energy Center (CEA) and Scientific and Technical Center for Buildings (CSTB). The authors also acknowledge the support of CNRS/INSIS (Cellule \'Energie) under the program ``\textit{Projets Exploratoires} --- 2017''.


\bigskip
\addcontentsline{toc}{section}{References}
\bibliographystyle{abbrv}
\bibliography{biblio}
\bigskip

\end{document}